\documentclass[a4paper,11pt,pdftex,bibtotocnumbered, headsepline,normalheadings, 
]
{article}
\usepackage{amssymb,amsmath,amsfonts,amsthm}
\usepackage[T1]{fontenc}
\usepackage{epsfig}
\usepackage{graphicx}
\usepackage{makeidx}

\usepackage{rotating}
\input xy 
\xyoption{all}
\xyoption{poly}
\usepackage{lscape}
\usepackage[colorlinks=true,linkcolor=black,pdftex]{hyperref}
\usepackage[hyperpageref]{backref} 
\usepackage{calc}
\theoremstyle{plain}

\newtheorem{satz}{Theorem}[section]

\newtheorem{coro}{Corollary}[section]
\newtheorem{lemma}{Lemma}
\newtheorem{prop}{Proposition}
\theoremstyle{definition}
\newtheorem{defi}{Definition}
\theoremstyle{definition}

\newtheorem*{rem}{Remark}

\theoremstyle{remark}

\DeclareMathOperator{\Hom}{\rm Hom}
\DeclareMathOperator{\End}{\rm End}

\DeclareMathOperator{\inj}{\hookrightarrow}

\newcommand{\oTo}{\xymatrix{ \ar@{^{(}->}[r]|{\mathbf{O}}& }}
\newcommand{\cTo}{\xymatrix{ \ar@{^{(}->}[r]|{\mathbf{|}}& }}
\newcommand{\coTo}{\xymatrix{ \ar@{^{(}->}[r]|{\mathbf{O}}|{\mathbf{|}}& }}

\DeclareMathOperator{\Bild}{Im}

\DeclareMathOperator{\Ind}{Ind}
\DeclareMathOperator{\Lie}{Lie}
\DeclareMathOperator{\eu}{eu}

\newcommand{\si}{\sigma}

\newcommand{\la}{\lambda}
\newcommand{\La}{\Lambda}
\newcommand{\al}{\alpha}

\newcommand{\De}{\Delta}

\newcommand{\B}{\mathbb{B}}
\newcommand{\C}{\mathbb{C}}

\newcommand{\G}{\mathbb{G}}

\newcommand{\N}{\mathbb{N}}

\renewcommand{\P}{\mathbb{P}}

\renewcommand{\S}{\mathbb{S}}

\newcommand{\U}{\mathbb{U}}

\newcommand{\W}{\mathbb{W}}

\newcommand{\Z}{\mathbb{Z}}

\newcommand{\mcB}{\mathcal{B}}
\newcommand{\mcC}{\mathcal{C}}

\newcommand{\mcE}{\mathcal{E}}

\newcommand{\mcG}{\mathcal{G}}
\newcommand{\mcH}{\mathcal{H}}

\newcommand{\mcP}{\mathcal{P}}

\newcommand{\mcR}{\mathcal{R}}

\newcommand{\mcT}{\mathcal{T}}
\newcommand{\mcU}{\mathcal{U}}

\newcommand{\mcZ}{\mathcal{Z}}

\begin{document}
\section*{Generalized quiver Hecke algebras }\label{Steinberg algebra}
\begin{center}Julia Sauter, University of Leeds,  \today \end{center}
\begin{abstract}
We generalize the methods of Varagnolo and Vasserot, \cite{VV} and partlially \cite{VV2}, 
to generalized quiver representations introduced by Derksen and Weyman in \cite{DW}. 
This means we have a general geometric construction of an interesting class of algebras (the Steinberg algebras for generalized quiver-graded Springer theory) containing skew group rings of Weyl groups with polynomial rings, (affine) nil Hecke algebras and KLR-algebras (=quiver Hecke algebras). Unfortunately this method works only in the \emph{Borel case}, i.e. all parabolic groups in the construction data of a Springer theory are Borel groups. Nevertheless, we try to treat also the parabolic case as far as this is possible here. 
\end{abstract}

This is a short reminder of Derksen and Weyman's generalized quiver representations from \cite{DW}.
\begin{defi}
A \emph{generalized quiver with dimension vector} is a triple $(\mathbb{G},G,V)$ where $\mathbb{G}$ is a reductive group, $G$ is a centralizer of a Zariski closed abelian reductive subgroup $H$ of $\mathbb{G}$,i.e. 
\[ G= C_{\G}(H) =\{g\in \G \mid ghg^{-1} =h\;\;\; \forall h\in H \} \]  
 (then $G$ is also reductive, see lemma below) and $V$ is a representation of $G$ which decomposes into irreducible representations which also appear 
in $\mathcal{G}:= \Lie (\mathbb{G})$ seen as an $G$-module.\\
A \emph{generalized quiver representation} is a quadruple $(\mathbb{G},G,V,Gv)$ where $(\mathbb{G},G,V)$ is a generalized quiver with dimension vector, $v$ in $V$ and $Gv$ is the $G$-orbit.  
\end{defi}

\begin{rem}
Any such reductive abelian group is of the form $H=A\times S$ with $A$ finite abelian and $S$ a torus, this implies that there exists finitely many elements $h_1,\ldots ,h_m$ such that $C_{\G} (H)=\bigcap_{i=1}^m C_{\G}(h_i)$, see for example Humphreys' book \cite{Hu}, Prop. in 16.4, p.107.
%
%
%
\end{rem}

We would like to work with the associated Coxeter systems, therefore it is sensible to assume $\G$ connected and replace $G$ by its identity component $G^o$. 
There is the following proposition

\begin{prop}
Let $\G$ be a connected reductive group and $H\subset \G$ an abelian group which lies in a maximal torus. We set $G:= C_G(H)^o=(\bigcap_{i=1}^m C_{\G}(h_i))^o$. Then it holds 
\begin{itemize}
\item[(1)] For any maximal torus $T\subset \G$, the following three conditions are equivalent:
\begin{itemize}
\item[(i)] $T\subset G$.
\item[(ii)] $H\subset T$.
\item[(iii)] $\{h_1,\ldots ,h_m\} \subset T$.
\end{itemize} 
\item[(2)] $G$ is a reductive group. 
\item[(3)] If $\underline{\Phi}$ is the set of roots of $\G$ with respect to a maximal torus $T$ with $H\subset T$, then $\Phi:= \{ \al \in \underline{\Phi} \mid \al (h) =1 \;\;\; \forall h\in H\} $ is the set of roots for $G$ with respect to $T$, its Weyl group is $\langle s_{\al}\mid \al \in \Phi\rangle $ and for all $\al\in \Phi$ the weight spaces are equal $\mathfrak{g}_{\al}=\mathcal{G}_{\al}$ (and $1$-dimensional $\C$-vector spaces). 
\item[(4)] There is a surjection 
\[ 
\begin{aligned}
\{ \B \subset \G \mid \B \text{ Borel subgroup, }H \subset \B \} &\to \{ B\subset G  \mid B \text{ Borel subgroup }\}\\
\B & \mapsto \B \cap G
\end{aligned}
\]
If $\underline{\Phi}^+$ is the set of positive roots with respect to $(\G, \B ,T)$ with $H\subset T$, then $\Phi^+:= \Phi\cap \underline{\Phi}^+$ is the set of positive roots for $(G, G\cap \B ,T)$ .
\end{itemize}
\end{prop}

\paragraph{proof:} Ad (1): This is easy to prove directly.\\
(2)-(4) are proven if $G=C_{\G}(h)^o$ for one semisimple element $h\in \G$ in Carters book \cite{C}, section 3.5.
 p.92-93. In general $G=(\bigcap_{i=1}^m C_{\G}(h_i)^o)^o$ for certain $h_i\in H, 1\leq i\leq m$. The result follows via induction on $m$. 
 Set $\G_1:= C_G(h_1)^o$. It holds $G=(\bigcap_{i=2}^m C_{\G_1}(h_i)^o)^o=C_{\G_1}(H)^o\subset \G_1$ and $\G_1$ is a connected reductive group. By induction hypothesis, all statements are true for $(G, \G_1)$, so in particular $G$ is a reductive group. The other statements are then obvious.
\hfill $\Box$

\subsubsection{Notational conventions}
We fix the ground field for all algebraic varieties and Lie algebras to be $\C$.\\
For a Lie algebra $\mathfrak{g}$ we define the $k$-th power inductively by $\mathfrak{g}^1:= \mathfrak{g}, \mathfrak{g}^k=[\mathfrak{g}, \mathfrak{g}^{k-1}]$. 
If we denote an algebraic group by double letters (or indexed double letters)  like $\mathbb{G}, \mathbb{B},\mathbb{U},...$ (or $\mathbb{G}^\prime$, $\mathbb{P}_J$, etc.)  we take the calligraphic letters for the Lie algebras, i.e. 
$\mcG, \mcB,\mcU, ...$ (or $\mcG^\prime, \mcP_J,$ etc) respectively. 
If we denote an algebraic group by roman letters (or indexed roman letters) like 
$G,B,U,...$ (or $G^\prime, P_J$, etc.) we take the small frakture letters for the Lie algebras, i.e. 
$\mathfrak{g}, \mathfrak{b},\mathfrak{u},.. $ (or $\mathfrak{g}^\prime, \mathfrak{p}_J$) respectively. \\
If we habe a subgroup $P\subset G$ of a group and an element $g\in G$ we write ${}^g\! P:= gPg^{-1}$ for the conjugate subgroup.

\noindent
We also recall the following. 
\begin{rem}
Let $(W,S)$ be a Coxeter system, $J\subset S$. Then $(W_J:= \langle J\rangle , J)$ is again a Coxeter system with the length function is the restriction of the length function of $(W,S)$ to elements in $W_J$. Then, the set $W^J$ of minimal length coset representatives $W^J\subset W$ for $W/W_J$ is defined via:  An element $w$ lies in $W^J$ if and only if for all $s\in J$ we have $l(ws)>l(w)$. Also there is a factorization $W=W^JW_J$ and if $w=xy$ with $x\in W^J,y\in W_J$, their lengths satisfy $l(w)=l(x)+l(y)$. We will fix the bijection $c_J\colon W^J\to W/W_J, w\mapsto wW_J$. The Bruhat order of $(W,S)$ can be restricted to $W^J$ and transferred via the bijection to $W/W_J$.\\
For two subsets $K,J\subset S$ define ${}^K W^J:= ({W^K})^{-1}\cap W^J$, the projection $W\to W_K \backslash W / W_J $ restricts to a bijection $ {}^KW^J\to W_K \backslash W / W_J $. \\
Let $(G, B,T)$ be a reductive group with Borel subgroup and maximal torus and $(W,S)$ be its associated Coxeter system. We fix for any element in $W$ a lift to the group $G$ and denote it by the same letter. 

\end{rem}

\subsection{Generalized quiver-graded Springer theory} 
We define a generalized quiver-graded Springer theory for generalized quiver representations in the sense of Derksen and Weymann. Given $(\mathbb{G}, \mathbb{P}_J, \mcU , H, V)  $ (and some not mentioned $H\subset T\subset \mathbb{B}\subset \mathbb{P}_J$)  with  
\begin{itemize}
\item[*] $\mathbb{G}$ is a connected reductive group, $H\subset T$ is a subgroup of a maximal torus in $\G$, we set 
$G=C_{\G}(H)^o$ (then $G$ is also reductive with $T\subset G$ is a maximal torus in $G$). 
\item[*] $T\subset \mathbb{B} \subset \mathbb{G}$ a Borel subgroup, then $B:= \mathbb{B}\cap G$ is a Borel subgroup of $G$,\\
         We write $(\mathbb{W},\mathbb{S})$ for the Coxeter system associated with $(\mathbb{G}, \mathbb{B}, T)$ and $(W,S)$ for the one associated 
         to $(G,B,T)$. Observe, that $W\subset \mathbb{W}$. For any $J\subset \mathbb{S}$ we set $\mathbb{P}_J:= \mathbb{B} \langle J\rangle \mathbb{B}$ and call it a  
         \textbf{standard parabolic group}.  
\item[*] Now fix a subset $J\subset \mathbb{S}$. We call a $\mathbb{P}_J$-subrepresentation $\mcU^\prime \subset  \mcG=\Lie(\G)$ (of the adjoint representation which we denote by $(g, x)\mapsto  {}^gx$, $g\in \G, x\in \mcG$) \emph{suitable} if  
\begin{itemize}
\item[$\bullet $] $(\mcU^{\prime})^T=\{0\}$,
\item[$\bullet $] $\mcU^{\prime}\cap {}^s\mcU^{\prime}$ is $\P_J$-stable for all $s\in \S$.
\end{itemize}
Let $\mcU=\bigoplus_{k=1}^r\mcU^{(k)}$ a $\mathbb{P}_J$-representation with each $\mcU^{(k)}$ is suitable. 
(Examples of suitable $\mathbb{P}_J$-representations are given by $\mcU^{\prime}=\mcU_{J^{\prime}}^{t}$ where $J\subset J^{\prime} \subset \S$, $\mcU_{J^{\prime}}=\Lie(\U_{J^{\prime}})$ with $\U_{J^{\prime }}\subset \P_{J^{\prime }}$ is the unipotent radical and 
$\mcU_{J^{\prime }}^{t}$ is the $t$-th power, $t\in \N$).          
We define $\mathbb{W}_{J} :=\langle J\rangle$ and $\mathbb{W}^J$ be the set of minimal coset representatives in $\mathbb{W}/\mathbb{W}_J$, 
$I_J:=W\backslash \mathbb{W}^J \subset W\backslash \mathbb{W}$ and 
\[
\bigsqcup_{J\subset \mathbb{S}} I_J 
\]
We call $I:=I_{\emptyset}$ the set of \textbf{complete dimension filtrations}. 
Let $\{ x_i\in \mathbb{W}\mid i\in I_J \}$ be a complete representing system of the cosets in $I_J$. 
Every element of the Weyl groups $\mathbb{W}$ (and $ W$) we lift to elements in $\mathbb{G}$ (and $G$) and denote the lifts by the same letter. 
For every $i\in I_J$ we set 
\[ P_{i} := {}^{x_i}\!\mathbb{P}_J \cap G, \]         
Observe that $H\subset T={}^{w}T \subset {}^w\P_J$ for all $w\in \W$, therefore ${}^w\P_J \cap G$ is a parabolic subgroup in $G$ for any $w\in \W$.
\item[*] $V=\bigoplus_{k=1}^r V^{(k)}$ with $V^{(k)} \subset \mcG$ is a $G$-subrepresentation.\\
$F_i=\bigoplus_{k=1}^r F^{(k)}_i$ with $F^{(k)}_i := V^{(k)} \cap {}^{x_i}\mcU^{(k)}$ is a $P_i$-subrepresentation of $V^{(k)}$. 
\end{itemize}
We define 
\[ 
\xymatrix{
& E_i := G\times^{P_i}F_i\ar[ld]_{\pi_i} \ar[rd]^{\mu_i}  &    &&  &\overline{(g,  f)} \ar[dl]\ar[dr]& \\
V && G/P_i && gf & & g P_i. 
}
\]
Now, there are closed embeddings $\iota_i\colon G/P_i\to \mathbb{G}/\mathbb{P}_J, \;  gP_i\mapsto   g x_i\mathbb{P}_J$ with for any $i\neq i^\prime $ in $I_J$ it holds $\Bild \iota_i \cap \Bild \iota_{i^\prime} =\emptyset$. Therefore, we can see $\bigsqcup_{i\in I_J} G_i/P_i$ as a closed subscheme of $\mathbb{G}/\mathbb{P}_J$. It can be identified with the closed subvariety of the fixpoints under the $H$-operation $(\mathbb{G}/\mathbb{P}_J)^H=\{ g\mathbb{P}_J\in \mathbb{G}/\mathbb{P}_J\mid hg\mathbb{P}_J= g\mathbb{P}_J \text{ for all }h\in H\}$.   
\[ 
\xymatrix{
& E_J:= \bigsqcup_{i\in I_J}E_i \ar[ld]_{\pi_J} \ar[rd]^{\mu_J}   \\
V && \mathbb{G}/\mathbb{P}_J 
}
\]
We also set 
\[ 
\xymatrix{
 &Z_{ij}:= E_i \times_{V} E_j \ar[dl]_{p_{ij}} \ar[dr]^{m_{ij}} &  &  &Z_J:= \bigsqcup_{i,j\in I_J} Z_{ij}\ar[dl]_{p_J} \ar[dr]^{m_J}&\\
 V && G/P_i\times G/P_j &  V&& \mathbb{G}/\mathbb{P}_J\times \mathbb{G}/\mathbb{P}_J. 
}
\]
In an obvious way all maps are $G$-equivariant. We are primarily interested in the following \textbf{Steinberg variety} 
\[ Z:= Z_{\emptyset} .\] 

The equivariant Borel-Moore homology of a Steinberg variety together with the convolution operation (defined by Ginzburg) defines a finite dimensional graded $\mathbb{C}$-algebra. We set   
\[ \mcZ_G:= H^G_*(Z) \]
which we call \textbf{(G-equivariant) Steinberg algebra}. The aim of this section is to describe $\mcZ_G$ in terms of generators and relation (for $J=\emptyset$). 
This means all $P_i$ are Borel subgroups of $G$. \\
If we set 
\[ H_{[p]}^G(Z) := \bigoplus_{i,j \in I} H_{e_i+e_j-p}^G(Z_{i,j}),\quad e_i =\dim_{\C} E_i
\]
then $H_{[*]}^G(Z)$ is a graded $H_G^*(pt)$-algebra
Then, we denote the right $\W$-operation on $I= W \setminus \W$ by $(i,w)\mapsto iw$, $i\in I, w\in \W$. We prove the following. 

\begin{satz}
Let $J=\emptyset$. Then $\mcZ_G\subset \End_{\C[\mathfrak{t}]^W-mod} (\bigoplus_{i\in I} \mcE_i), \;\mcE_i=\C[\mathfrak{t}]= \C[x_i(1), \ldots , x_i(n)], i\in I$ is the $\C $-subalgebra generated by 
\[1_i, i\in I, \quad  z_i(t), 1\leq t\leq n=rk(T), i\in I, \quad  \si_i(s), s\in \S, i\in I\] 
defined as follows for $k\in I$, $f\in \mcE_k$.
\[
\begin{aligned}
{1_i} (f) &:= \begin{cases} f,  &\text{ if }i=k, \\
                                0, &\text{ else. }
\end{cases} \\
{z_{i}(t)}(f)&:= \begin{cases} x_i(t)f ,  &\text{ if }i=k, \\
                       0, & \text{ else. }
\end{cases} \\
{\sigma_{i} (s)} (f) &:= \begin{cases} q_i(s)\frac{s(f)-f}{\al_s}   , \quad  &(\in  \mcE_i)\text{ if }i=is=k, \\
                                                 q_{i}(s) s(f)                  \quad  & (\in \mcE_i) \text{ if }i\neq is=k, \\
                       0, & \text{ else. }
\end{cases}                     
\end{aligned}
\] 
where
\[q_i(s):= \prod_{\al \in \Phi_{\mcU}, s(\al) \notin \Phi_{\mcU}, x_i(\al) \in \Phi_V} \al \quad \quad  \in \mcE_i.\]
and $\Phi_{\mcU} =\bigsqcup_{k} \Phi_{\mcU^{(k)}}$, $\Phi_{\mcU^{(k)}} \subset \Hom_{\C}(\mathfrak{t}, \C)\subset \C[\mathfrak{t}]$ is the set of $T$-weights 
for $\mcU^{(k)}$ and $\Phi_V=\bigsqcup_{k} \Phi_{V^{(k)}}$, $\Phi_{V^{(k)}}\subset \Hom_{\C}(\mathfrak{t}, \C)$ is the set of $T$-weights for $V^{(k)}$. \\
Furthermore, it holds 
\[\deg 1_i =0, \; \deg z_i(k) =2, \; \deg \si_i(s) = \begin{cases} 2(\deg q_i(s))-2, & \text{ if } is=i \\
                                                                  2\deg q_i(s) , &\text{ if }is\neq i \end{cases}      \]
where $\deg q_i(s)$ refers to the degree as homogeneous polynomial in $\C[\mathfrak{t}]$.      
\end{satz}

The generality of the choice of the $\mcU$ in the previous theorem is later used to understand the case of an arbitrary $J$ as a an algebra of the form 
$e_J\mcZ_G e_J$ for an associated Borel-case Steinberg algebra $\mcZ_G$ and $e_J$ an idempotent element (this is content of a later article called \emph{parabolic Steinberg algebras}).\\ 
For $J=\emptyset, \mcU=\Lie (\U)^{\oplus r}$ for $\U\subset \B$ the unipotent radical 
we have the following result which generalizes KLR-algebras to arbitrary connected reductive groups and allowing quivers with loops.  

\begin{coro}
Let $J=\emptyset, \mcU=Lie (\U)^{\oplus r}$, $\U\subset \B$ the unipotent radical. Then \[\mcZ_G\subset \End_{\C[\mathfrak{t}]^W-mod} (\bigoplus_{i\in I} \mcE_i),\]
$\mcE_i=\C[\mathfrak{t}]= \C[x_i(1), \ldots , x_i(n)], i\in I$ is the $\C$-subalgebra generated by 
\[1_i, i\in I, \quad  z_i(t), 1\leq t\leq n=rk(T), i\in I, \quad  \si_i(s), s\in \S, i\in I.\] 
Let $f\in \mcE_k, k\in I$, $\;\al_s \in \Phi^+$ be the positive root such that $s(\al_s)=-\al_s$. It holds 
\[
\begin{aligned}
{\sigma_{i} (s)} (f) &:= \begin{cases} \al_s^{h_i(s)}\frac{s(f)-f}{\al_s} , \quad  &\text{ if }i=is=k, \\
                                                    \al_s^{h_{i}(s)} s(f) &\text{ if }i\neq is =k, \\
                       0, & \text{ else. }
\end{cases}                     
\end{aligned}
\]
where 
\[ h_i (s):= \# \{ k\in \{1,\ldots ,r\} \mid x_i(\al_s)\in \Phi_{V^{(k)}}\}\]
where $V=\bigoplus_k V^{(k)}$ and $\Phi_{V^{(k)}}\subset \underline{\Phi}$ are the $T$-weights of $V^{(k)}$.  
\begin{itemize}
\item[(1)] If $Wx_i\neq Wx_is$ then 
\[h_i(s)=\# \{ k \mid V^{(k)} \subset \mcR, x_i(\al_s)\in \Phi_{V^{(k)}}\}. \]
We say that this number \emph{counts arrows}. 
\item[(2)] If $Wx_i=Wx_is$, then 
 \[h_i(s)=\# \{ k \mid V^{(k)} \subset \mathfrak{g}, x_i(\al_s)\in \Phi_{V^{(k)}}\}. \]
We say that this number \emph{counts loops}.
\end{itemize}
\end{coro}

In the case of the previous corollary we call the Steinberg algebra $\mcZ_G$ \textbf{generalized quiver Hecke algebra}. 
It can be described by the following generators and relations.  For a reduced expression $w=s_1s_2\cdots s_k$ we set 
\[
\si_i (s_1s_2\cdots s_k) :=
\si_{i} (s_1) \si_{is_1}(s_2)\cdots \si_{is_1s_2\cdots s_{k-1}} (s_{k})
\]
Sometimes, if it is understood that the definition depends on a particular choice of a reduced expression for $w$, we write $\si_i (w):= \si_i (s_1s_2\cdots s_k)$. 
Furthermore, we consider \[\Phi\colon \bigoplus_{i\in I} \C[x_i(1),\ldots x_i(n)]\cong \bigoplus_{i\in I} \C[z_i(1),\ldots z_i(n)], \; x_i(t)\mapsto z_i(t)\] 
as the left $\W$-module $\Ind_W^{\W} \C[\mathfrak{t}]$, we fix the polynomials 
\[ c_{i}(s,t):= \Phi(\si_i(s)(x_i(t))) \;\in \bigoplus_{i\in I} \C[z_i(1),\ldots z_i(n)], \quad i\in I, \; 1\leq t\leq n,\; s\in \S. \]

Now, we can describe under some extra conditions the relations of the generalized quiver Hecke algebras. 
\begin{prop} Under the following assumption for the data $(\G, \B , \mcU= (\Lie (\U))^{\oplus r}, H,V), J=\emptyset\colon $
Let $\S\subset \W=Weyl(\G,T)$ be the simple reflections, we assume for any $s,t\in \S$ 
\begin{itemize}
\item[(B2)] If the root system spanned by $\al_s, \al_t$ is of type $B_2$ (or $stst=tsts$ is the minimal relation), then for every $i\in I$ such that 
$is=i=it$ it holds $h_i(s), h_i(t)\in \{0,1,2\}$. 
\item[(G2)] If the root system spanned by $\al_s, \al_t$ is of type $G_2$ (or $ststst=tststs$ is the minimal relation), then for every $i\in I$ such that 
$is=i=it$ it holds $h_i(s)=0=h_i(t)$. 
\end{itemize}  
Then the generalized quiver Hecke algebra for $(\G, \B , \mcU= (Lie (\U))^{\oplus r}, H,V), J=\emptyset$ is the $\C$-algebra with generators 
 \[1_i, i\in I, \quad  z_i(t), 1\leq t\leq n=rk(T), i\in I, \quad  \si_i(s), s\in \S, i\in I\] 
 and relations 

\begin{itemize}

\item[(1)] (\emph{orthogonal idempotents})
\[
\begin{aligned}
1_i 1_j  &=\delta_{i,j} 1_i, \\
1_i z_i(t) 1_i  &= z_i(t), \\
1_i \si_i(s) 1_{is} &= \si_i(s)
\end{aligned}
\]
\item[(2)] (\emph{polynomial subalgebras}) 
\[z_i(t) z_i(t^\prime)= z_i(t^\prime) z_i(t)\]

\item[(3)] (\emph{ relation implied by $s^2=1$}) 
\[\si_i(s) \si_{is}(s) = \begin{cases}  0 &, \text{ if }is=i,\;  h_i(s) \text{ is even } \\
                                    -2\al_s^{h_i(s)-1} \si_i(s) & , \text{ if }is=i,\;  h_i(s) \text{ is odd } \\
                                    (-1)^{h_{is}(s)} \al_s^{h_i(s)+h_{is}(s)} &, \text{ if } is\neq i
\end{cases}
\]

\item[(4)] (\emph{straightening rule})\\
\[ 
\si_i(s) z_i(t) - s(z_i(t))\si_i(s) = \begin{cases} c_{i}(s,t), &\text{, if } is=i \\
0 &\text{, if } is\neq i. \end{cases}
\]

\item[(5)] (\emph{braid relations})\\
Let $s,t\in \S, st=ts$, then 
\[  \si_i(s)\si_{is}(t) =  \si_i (t) \si_{it}(s)
\]
 Let $s, t\in \S$ not commuting such that $x:=sts \cdots= tst\cdots $ minimally, $i\in I$. 
There exists explicit polynomials $(Q_w)_{w< x}$ in $\al_s, \al_t\in \C[\mathfrak{t}]$ such that 
\[ 
\si_i (sts \cdots ) -\si_i (tst\cdots ) = \sum_{w<x} Q_w \si_i (w) 
\]
(observe that for $w<x$ there exists just one reduced expression).  

\end{itemize}

\end{prop}

The proof you find in the end, see Prop. \ref{relations}.

\subsubsection{Relationship between parabolic groups in $G$ and $\G$} \label{bigSmall}

For later on, we need to understand the relationship between parabolic subgroups in $\G$ and in $G$. 
Recall that a parabolic subgroup is a subgroup which contains a Borel subgroup, every parabolic subgroups is conjugated to a standard parabolic subgroup. The standard parabolic subgroups wrt $(G,B,T)$ are in bijection with the set of subsets of $S$, via $J\mapsto B\langle J\rangle B=: P_J$. As a first step, we need to study the relationship of the Coxeter systems $(W,S)$ and $(\W, \S)$.

\begin{lemma}
It holds $G\cap \W = W$. It holds $W\cap \S \subset S$. Let $l_S$ be the length function with respect ot $(W,S)$ and $l_{\S}$ be the length function with respect to $(\W, \S)$. For every $w\in W$ it holds  
$l_{S}(w)\leq l_{\S}(w)$.
\end{lemma}

\paragraph{proof:}
$N_{\G} (T)\cap G = N_G(T)$ implies $G\cap \W =W$. 
The inculsion $ \Phi^+ \cap s (-\Phi^+) \subset  \underline{\Phi}^+ \cap s (-\underline{\Phi}^+)$ for any $s\in \S$ implies $W\cap \S \subset S$. \\
Let $w= t_1\cdots t_r\in W$, $t_i\in S$ reduced expression and assume $l_{\S}(w) <r$. It must be possible in $\W$ to write $w$ as a subword of 
$t_1\cdots \hat{t_i}\cdots t_r$ for some $i\in \{1,\ldots , r\}$. But then $r= l_{S}(w)\leq l_S( t_1\cdots \hat{t_i}\cdots t_r) <r$. 
\hfill $\Box$

\begin{defi}
We call $J\subset \S$. We say that $J$ is \textbf{$S$-adapted} if for all $s\in S$ with $s=s_1\cdots s_r$ a reduced expression in $(\W ,\S)$ such that there exists $i\in \{ 1,\ldots , r\}$ with $s_i\in J$ then it also holds $\{ s_1,\ldots , s_r\} \subset J$. 
\end{defi}
%
%
 
\begin{lemma} \label{cap} 
\begin{itemize}
\item[(a)] 
Intersection with $G$ defines a map 
\[
\begin{aligned}
 \{ \P_J \mid J\subset \S \text{ is }S-\text{adapted} \} &\to \{ P_J \mid J\subset S\} \\
 \P_J &\mapsto \P_J \cap G = P_{S\cap\W_J} 
 \end{aligned}
 \]
\item[(b)] 
Let $G\cap{}^x\B$ is a Borel subgroup of $G$ with $\B\subset G$ a Borel subgroup and $x\in \W$. Let $s\in \S$, then it holds 
\begin{itemize}
\item[(1)] If $Wxs \neq Wx$ then $G\cap {}^{xs}\B =G\cap {}^x\B$.
\item[(2)] If $Wxs =Wx$, then ${}^x s\in W$ and $G\cap {}^{xs}\B= {}^{{}^x s} [ G\cap {}^x\B] $. 
\end{itemize}
This gives an algorithm to find for any $x\in \W$ a $z\in W$ such that $G\cap {}^x\B = {}^z[G\cap \B]$. \\
Also, for every $J\subset \mathbb{S}$ it then holds $G\cap {}^x\P_J  = {}^z[G\cap\P_J]$  and $W\cap {}^x\W_J = {}^z [W\cap \W_J]$ where $x\in \W, z\in W$ as before and for every $S$-adapted $J\subset \mathbb{S}$ 
\[ G\cap {}^x \P_J = {}^z \P_{S\cap \W_J}. \]  
\end{itemize}
\end{lemma} 

\paragraph{proof:} 
\begin{itemize}
\item[(a)]
It holds by the previous lemma $G\cap \W_J = W\cap \W_J$ and because $J$ is $S$-adapted it holds $ W\cap \W_J= \langle S\cap \W_J \rangle$, to see that:\\
 Let $w=t_1\cdots t_r\in \W_J$ with $t_i\in S$ an $S$-reduced expression, we need to see $t_i\in \W_J, 1\leq i\leq r$. Wlog assume $t_1\notin \W_J$. As $J$ is $S$-adapted, there exists a $\S$-reduced expression with elements in $J$ of $w$ which is a subword of $t_2\cdots t_r$. But this means a word of $S$-length $r$ is a subword of a word of $S$-length $r-1$, therefore $t_1\in \W_J$.\\
Now, the following inclusion is obvious 
\[ P_{S\cap \W_J} = B \langle G\cap \W_J \rangle B \subset G\cap \P_J.\]
Because $B\subset \P_J \cap G$ there has to exist $(\W_J\cap S) \subset J^\prime \subset S$ such that $\P_J \cap G=P_{J^\prime}$, we need to see $(S\cap \W_J)=J^\prime$.
Let $s\in J^\prime $, then $s\in \P_J=\B \W_J \B$ implies $s\in \W_J$. \\
\item[(b)] 
Let $s\in \S, {}^x s\notin W$, then $\pm x(\al_s)\notin \Phi$ and this implies
\[ \Phi \cap xs(\underline{\Phi})= \Phi \cap [x(\underline{\Phi})\setminus \{ x(\al_s)\} \cup \{ -x(\al_s)\}]=\Phi \cap x(\underline{\Phi}).\]
Therefore, the Lie algebras of the Borel groups $G\cap{}^x\B$ and $G\cap{}^{xs}\B$ have the same weights for $T$, this proves they are equal. \\
The point (2) is obvious. 
\end{itemize}
\hfill $\Box $

\begin{rem} \label{minis} 
In the setup of the beginning, we can always find unique representatives $x_i\in \W, i\in I$ for the elements in $W\setminus\W$ which fulfill 
\[ 
B_i=G\cap {}^{x_i}\B = G\cap \B =B.
\]
This follows because for every $i\in I$ there is a bijection 
\[
\begin{aligned}
Wx_i &\to \{ \text{ Borel subgroups of }G \text{ containing }T\}\\
vx_i &\mapsto {}^v[G\cap {}^{x_i}\B]
\end{aligned}
\]
Then, there exists a unique $v\in W$ such that ${}^v[G\cap {}^{x_i}\B]=G\cap \B$, replace $x_i$ by $vx_i$ as a representative for $Wx_i$.\\
We will call these representatives \textbf{minimal coset representatives}\footnote{if $G$ is a Levi-group in $\G$ they are the minimal coset representatives, in this more general situation the notion is not defined.}. Observe for $is\neq i$ it holds $x_{is}=x_is$ by lemma \ref{cap}, (b), (2). \\
But since the images of $G/B_i, i\in I$ inside $\G/\B$ are disjoint, we prefer not to identify all $B_i, i\in I$. \\
In general, in the parabolic setup, it holds $P_i\neq P_j$ for $i\neq j$. 
\end{rem}

\begin{lemma} (factorization lemma)
Let $J,K\subset \S$ be $S$-adapted and set $L:=S\cap \W_J,M:=S\cap \W_K$.  
\begin{itemize} 
\item[(1)] It holds $W^L = W\cap \W^J$ and for every element in $w\in W$ the unique decomposition as $ w= w^J w_J$, $w^J  \in \W^J, w_J \in
\W_J$ 
fulfills 
$w^J\in  W^L = W\cap  \W_J$ , $w_J\in  W_L=W\cap \W_J$.
\item[(2)] It holds ${}^J\W^K \cap W = {}^L W^{M}$. In particular, every double coset $\W_J w\W_K$ with $w\in W$ contains a unique element of ${}^{L} W^{M}$.
\end{itemize}
\end{lemma}

\paragraph{proof:}
\begin{itemize}
\item[(1)] It holds $W_L (W\cap \W_J) =W =W\cap \W^J\W_J \supset (W\cap \W^J)(W\cap \W^J)$, the uniqueness of the factorization in $W$ implies $(W\cap \W^J)\subset W^L$.\\
Now take $a\in W^L$, we can factorize it in $\W$ as $a=a^Ja_J$ with $a^J\in \W^J, a_J\in \W_J$. We show that $a_J\in W$. 
Write $a=t_1\cdots t_r$ $S$-reduced expression, assume $a_J\neq e$, 
then there exists a unique $i\in \{1,\ldots , r\}$ such that $a_J$ is a subword of $t_i\cdots t_r$ but no subword of $t_{i+1}\cdots t_r$. 
Then, $t_i$ must have a subword contained in $\W_J$, as $J$ is $S$-adapted we get $t_i\in \W_J$. Continue with $t_i^{-1}a_J$ being a subword of $t_{i+1}\cdots t_r$. By iteration you find $a_J= t_{i_1}\cdots t_{i_k}\in W$ for certain $i=i_1<\cdots <i_k$, $i_j\in \{1,\ldots , r\}$. This implies $a_J=e$ and 
$a=a^J\in W\cap \W^J$. 
\item[(2)] By definition ${}^J\W^K \cap W = (\W^J)^{-1}\cap \W^K \cap W = (W^L)^{-1} \cap W^M = {}^LW^M$. 
\end{itemize}
\hfill $\Box $

\subsubsection{The equivariant cohomology of flag varieties}

\begin{lemma} (\textbf{The (co)-homology rings of a point}) \\
Let $G$ be reductive group, $T\subset P\subset G$ with $P$ a parabolic subgroup and $T$ a maximal torus, we write $W$ for the Weyl group associated to $(G,T)$ and $ X(T)=\Hom_{Gr}(T, \C^*)$ for the group of characters. Let $ET$ be a contractible topological space with a free $T$-operation from the right.  
\begin{itemize}
\item[(1)] 
For every character $\la \in X(T)$ denote by 
\[ S_{\la} := ET\times^T \C_{\la} \] 
the associated $T$-equivariant line bundle over $BT:=ET/T$ to the $T$-representation $\C_{\la}$ which is $\C$ with the operation $t\cdot c := \la (t) c$. 
The first chern class defines a homomorphism of abelian groups 
\[ 
c\colon X(T) \to H^2(BT), \quad \la  \mapsto c_1 ( S_{\la}).
\] 
Let ${\rm{Sym}}_{\C}(X(T))$ be the symmetric algebra with complex coefficients generated by $X(T)$, it can be identified with the ring of regular function $\C[\mathfrak{t}]$ on $\mathfrak{t}=Lie(T)$ (with doubled degrees), where $X(T)\otimes_{\Z}\C$ is mapped via taking the differential (of elements in $X(T)$) to $\mathfrak{t}^*=\Hom_{\C-lin}(\mathfrak{t}, \C) \subset \C[\mathfrak{t}]$ (both are the degree 2 elements).   \\
  
The previous map extends to an isomorphism of graded $\C$-algebras
\[
 \C[ \mathfrak{t}] \to H^*_T (pt)=H^*(BT)\\
\] 

In fact this is a $W$-linear isomorphism where the $W$-operation on $\C[\mathfrak{t}]$ is given by, $(w,f)\mapsto w(f), w\in W, f\in \C[\mathfrak{t}]$ with 
\[ w(f) \colon \mathfrak{t}\to \C, t\mapsto f(w^{-1}tw).\] 
We can choose $ET$ such that it also has a free $G$-operation from the right (i.e. $ET:=EG$), then $BT=ET/T$ has an induced Weyl group action from the right given by 
$xT\cdot w:= xwT$, $w\in W, x\in ET$. The pullbacks of this group operation induce a left $W$-operation on $H_T^*(pt)$.  
 
\item[(2)] $H_T^*(pt) =H^T_{-*}(pt)$, $H_G^*(pt)= (H_T^*(pt))^W = (H^T_{-*}(pt))^W= H^G_{-*}(pt)$. 
\end{itemize}
\end{lemma} 

\paragraph{proof:}
\begin{itemize}
\item[(1)]
For the isomorphism see for example and the explanation of the $W$-operation see (L. Tu; Characteristic numbers of a homogeneous space, axiv, \cite{Tu})
\item[(2)] Use the definition and Poincare duality for the first isomorphism, for the second also use the splitting principle.

\end{itemize}

\hfill $\Box$
 
\begin{lemma} \label{CohomOfE}(\textbf{The cohomology rings of homogeneous vector bundles over $G/P$}) \\
Let $G$ be reductive group, $T\subset B\subset P\subset G$ with $B$ a Borel subgroup, 
$P$ parabolic and $T$ a maximal torus. 
\begin{itemize} 
\item[(1)] 
For $\la \in X(T)$ we denote be $L_{\la}:= G\times^B \C_{\la}$ 
the associated line bundle to the $B$-representation $\C_{\la}$ given by the trivial representation when restricted to the unipotent radical and 
$\la$ when restricted to $T$. Let $\mu\colon E\to G/B$ be a $G$-equivariant vector bundle.
Then, $\mu^*(L_{\la})$ is a line bundle on $E$ and  
\[ K_{\la}:= EG\times^G \mu^*(L_{\la}) \to EG\times^G E \]
is a line bundle over $EG\times^G E$. 
There is an isomorphism of graded $\C$-algebras
\[
\begin{aligned}
 \C[ \mathfrak{t}] &\to   H_G^*(E)=H^*(EG\times^G E) \\
  X(T)\ni \quad \la  & \mapsto c_1(K_{\la}).
 \end{aligned} 
 \] 
with $deg \la =2$ for $\la \in X(T)$. \\
(By definition, equivariant chern classes are defined as $c_1^G(\mu^* L_{\la}):= c_1(K_{\la})$).  
\item[(2)] 
Let $\mu\colon E\to G/P$ be a $G$-equivariant vector bundle, then there is an isomorphism of graded $\C$-algebras 
\[ H_G^*(E) \to (H_T^*(pt))^{W_L}. \]
\end{itemize}
\end{lemma}

\paragraph{proof:}
\begin{itemize}
\item[(1)] Arabia proved that $H_G^*(G/B)\cong H_T^*(pt)$ as graded $\C$-algebras (cp. \cite{Ar2}), the composition with the isomorphism from the previous lemma 
gives an isomorphism 
\[ 
c\colon \C[\mathfrak{t}] \to H_G^*(G/B), \colon \la\mapsto c_1(EG\times^G L_{\la})=:c_1^G(L_{\la})
\]
Now, we show that for a vector bundle $\mu\colon E\to G/P$ with $P\subset G$ parabolic, the induced pullback map 
\[ 
\mu^*\colon H_G^*(G/P) \to H_G^*(E),  \quad c_1^G(L_{\la})\mapsto c_1^G(\mu^*L_{\la})
\]
is an isomorphism of graded $H_G^*(pt)$-algebras. We already know that it is a morphism of graded $H_G^*(pt)$-algebras, to see it is an isomorphism, apply the definition and Poincare duality to get a commutative diagram 
\[
\xymatrix{
H_G^k(G/P) \ar[r]^{\mu^*}\ar[d]_{\cong} & H_G^k(E)\ar[d]_{\cong} \\
H^G_{2\dim G/P-k}(G/P)\ar[r]^{\mu^*} & H^G_{2\dim E -k}(E)
}
\]
the lower morphism $\mu^*$ is the pullback morphism which gives the Thom isomorphism, therefore the upper $\mu^*$ is also an isomorphism.  
\item[(2)]
By the last proof, we already know $H_G^*(E)\cong H_G^*(G/P)$. Then apply the isomorphism of Arabia see \cite{Ar2}, this gives $H_G^*(G/P)\cong H_P^*(pt)$. Now, 
$P$ homotopy-retracts on its Levy subgroup $L$, this implies $H_P^*(pt)=H_L^*(pt)$, together with the (2) in the previous lemma we are done. 
\end{itemize}
\hfill $\Box $

\begin{lemma} (\textbf{The cohomology ring of the flag variety as subalgebra of the Steinberg algebra})\\
 Let $G$ be reductive group, $T\subset P\subset G$ with $P$ parabolic and $T$ a maximal torus.  
Let $V$ be a $G$-representation and $F\subset V$ be a $P$-subrepresentation, let $E:= G\times^PF$ and $Z:= E\times_V E$ be the associated Steinberg variety.   
The diagonal morphism $E\to E\times E$ factorizes over $Z$ and induces an isomorphism $E\to Z_e$ which induces an isomorphism of algebras 
\[ H_G^*(G/P) \to H^G_{2\dim E -*}(Z_e), \] 
recall that the convolution product on $H_*^G(Z_e)$ maps degrees $(i,j)\mapsto i+j - 2\dim E $. 

\end{lemma}

\paragraph{proof:}
Obviously you have an isomorphism $H_G^*(G/B) \xrightarrow{\mu^*} H_G^*(E)\cong H_G^*(Z_e)\to H_{2\dim E -*}^G(Z_e)$ where the last isomorphism is Poincare duality. But we need to see that this is a morphism of algebras where $H_*^G(Z_e)$ is the convolution algebra with respect to the embedding $Z_e\cong E\xrightarrow{diag} E\times E$. This follows from \cite{CG}, Example 2.7.10 and section 2.6.15.

\hfill $\Box$

We observe that the algebra $\C[\mathfrak{t}]$ with generators $t\in \mathfrak{t}^*$ in degree $2$ plays three different roles in the last lemmata. 
It is the $T$-equivariant cohomology of a point, it is the $G$-equivariant cohomology of a complete flag variety $G/B$, it can be found as the subalgebra 
$H_*^G(Z_e)\subset H_*^G(Z)$.

\subsubsection{Computation of fixed points}
Recall the following result, for example see  \cite{Hae}, satz 2.12, page 13.
\begin{lemma}
Let $T\subset P\subset G$ be reductive group with a parabolic subgroup $P$ and a maximal torus $T$. Let $W$ be the Weyl group associated to $(G,T)$ and 
$Stab(P):= \{ w\in W \mid wPw^{-1}=P \}$. For $w= x Stab(P)\in W/  Stab (P) $ we set $wP:= xP \in G/P$. 
Then, it holds 
\[ (G/P)^T= \{ wP \in G/P \mid w \in W/ Stab (P) \} \] 
\end{lemma}

\begin{lemma} Let $P_1,P_2\subset G$ be a reductive group with two parabolic subgroup, $F_1, F_2\subset V$ a $G$-representation 
with a $P_1$ and $P_2$-subrepresentation. Assume $(GF_i)^T=\{0\}$. 
We write \\
$(E_i=G\times^{P_i}F_i, \mu_i \colon E_i\to G/P_i, \pi \colon E_i\to V)$ for the associated Springer triple and  \\
$Z:= E_1\times_V E_2, m\colon Z\to (G/P_1)\times (G/P_2)$ for the Steinberg variety. \\
Then, there are induced a bijections $\mu_i^T\colon E_i^T \to  (G/P_i)^T, m^T \colon Z^T \to (G/P_1)^T\times (G/P_2)^T$.
More explicit we have 
\[ 
\begin{aligned}
E_i^T &=\{\phi_w:= (0, wP_i) \in V\times G/P_i \mid w\in W/ Stab (P_i) \} \subset E_i \\
Z^T &=\{ \phi_{x,y} :=(0, xP_1, yP_2) \in V\times G/P_1\times G/P_2 \mid x\in W/ Stab (P_1),  y\in W/ Stab (P_2) \} \subset Z. 
\end{aligned}
\]  
Furthermore, for any $w\in W/ Stab (P_2)$ let $Z^w:= m^{-1}(G\cdot (P_1, wP_2))$ and 
$m_w:= m|_{Z^w}\colon Z^w \to G\cdot (P_1, wP_2)$ the induced map. There is an induced 
Bruhat order $\leq $ on $W/ Stab (P_2)$ by taking the Bruhat order of minimal length representatives. 
\[ 
\begin{aligned}
(Z^w)^T &= \{ \phi_{x, xw}= (0, xP_1, xwP_2) \in V\times G/P_1\times G/P_2 \mid x\in W\} \\
{\overline{Z^w}}^T &= \{ \phi_{x, xv} \mid x\in W, v\leq w \} = \bigcup_{v\leq w} (Z^v)^T 
\end{aligned}
\]
There is a bijection $W/(Stab (P_1)\cap {}^w Stab(P_2)) \to (Z^w)^T, \; x\mapsto \phi_{x,xw}$. 
\end{lemma}
 
 \paragraph{proof}
 Obviously, it holds $E_i^T \subset V^T \times (G/P_i)^T= \{ 0\} \times (G/P_i)^T$. But we also have a zero section 
 $s$ of the vector bundle $\pi \colon E_i\to G/P_i$ 
 which gives the closed embedding $G/P_i \to E_i \subset V\times (G/P_i), gP_i\mapsto (0,gP_i)$. \\
 It holds $Z^T \subset V^T \times (G/P_1)^T \times (G/P_2)^T =\{ 0\} \times (G/P_1)^T \times (G/P_2)^T $.
 But using the description of $Z=\{ (v, gP_1, hP_2 )\in V\times G/P_1\times G/P_2 \mid (v, gP_1)\in E_1, (v, hP_2) \in E_2\}$, 
 we see that $\{0\}\times (G/P_1)^T\times (G/P_2)^T \subset Z$ and these are obviously $T$-fixed points. \\
 We have $(Z^w)^T \subset Z^w\cap Z^T =\{ \phi_{x,xw}\mid x\in W\} $ and one can see the other inclusion, too. Also, we have $\overline{Z^w}^T \subset 
 (\bigcup_{v\leq w} Z^v)^T = \bigcup_{v\leq w} (Z^v)^T$. Consider the closed embedding 
 \[ s\colon G/P_1\times G/P_2 \to Z , \quad (gP_1, hP_2)\mapsto (0,gP_1, hP_2) .\]
Clearly $s(G(P_1,wP_2))\subset Z^w\subset \overline{Z^w}$, but since $s$ is a closed embedding we have 
\[ \bigcup_{v\leq w} (Z^v)^T\subset s(\overline{G(P_1,wP_2)})=\overline{s(G(P_1,wP_2))} \subset \overline{Z^w} \] 
which yields the other inclusion.              \hfill $\Box $

\paragraph{Notation for the fixed points}
Now, in the set-up of the beginning this gives the following:\\
 Observe, that $(\bigsqcup_{i\in I_J} G_i/P_i)^T =((\mathbb{G}/\mathbb{P}_J)^H)^T= (\mathbb{G}/\mathbb{P}_J)^T,$
and  $(\mathbb{G}/\mathbb{P}_J)^T = \{ w \mathbb{P}_J \mid w\in \mathbb{W}^J \}$. For any $w\in \W^J$ there exists a 
unique $i\in I_J$ such that $x^i:= wx_i^{-1}\in W$, this implies $w\P_J = x^i (x_i\P_J) = \iota_i(x^i P_i) \in (G/P_i)^T$. Therefore, we write 
\[
\begin{aligned}
(\bigsqcup_{i\in I_J} G_i/P_i)^T &= (\mathbb{G}/\mathbb{P}_J)^T =
\bigsqcup_{i\in I_J} \{ wx_i \mathbb{P}_J \mid w\in W/{}^{x_i}\mathbb{W}_J \subset \mathbb{W}/{}^{x_i}\mathbb{W}_J \}\\
E_J^T &=\{ \phi_{wx_i} = (0,wx_i \mathbb{P}_J )\mid i\in I_J, w\in W/{}^{x_i}\mathbb{W}_J\}   \\
Z_J^T &= \bigsqcup_{i,j\in I_J}\{ \phi_{wx_i, vx_j} =(0,wx_i \mathbb{P}_J , vx_j\mathbb{P}_J)\mid   w\in W/{}^{x_i}\mathbb{W}_J, v\in W/{}^{x_j}\mathbb{W}_J\}.
\end{aligned}
\]
Let $w,v\in \W^J$ and $i,j\in I_J$ such that $w^i:= wx_i^{-1}\in W, v^j:= vx_j^{-1} \in W$. 
We then set $\phi_w:= \phi_{w^i x_i}$, $\phi_{w,v}=\phi_{w^i x_i, v^j,x_j}$.\\ 
As we have bijections $(0,wx_i \mathbb{P}_J )\mapsto wx_i \mathbb{P}_J, (0,wx_i \mathbb{P}_J , vx_j\mathbb{P}_J)\mapsto (wx_i \mathbb{P}_J , vx_j\mathbb{P}_J) $ between   $E^T$ and $(\mathbb{G}/\mathbb{P}_J)^T$, $Z^T$ and $(\mathbb{G}/\mathbb{P}_J\times \mathbb{G}/\mathbb{P}_J)^T$, we denote the $T$-fixed by the same symbols. 


\subsubsection{The fibres over the fixpoints}

Remember, by definition we have $F_i=\mu_J^{-1}(\phi_{x_i})$. For any $w=w ^ix_i\in \W ^J, w^i\in W$ 
We set 
\[
F_w :=\mu_J^{-1}(\phi_{w})= 
\mu_i^{-1}(\phi_{w^{i}x_i}) = {}^{w^i}F_i  = \bigoplus_{k=1}^r V^{(k)}\cap {}^{w}\mcU^{(k)}
\]  
and if also $x \in \W / (\W_J \cap {}^w \W_J)$ (i.e. the definition does not depend on the choice of a representative in the coset)
\[
\begin{aligned}
F_{x,xw} &:= m_J^{-1} (\phi_{x,xw}) = 
F_x\cap F_{xw} \\
&=\bigoplus_{k=1}^r V^{(k)}\cap {}^{x}[\mcU^{(k)}\cap {}^{w}\mcU^{(k)}] 
\end{aligned}
\]
\paragraph{For $J=\emptyset, \mcU=\Lie(\U)^{\oplus r}\colon$} 
We choose $V=\bigoplus_{k=1}^t V^{(k)} \oplus \bigoplus_{k=t+1}^rV^{(k)}$ with 
$V^{(k)} \subset \mcR, 1\leq k\leq t$, $V^{(k)} =\mathfrak{g}^{(k)}$ with $\mathfrak{g}^{(k)}\subset \mathfrak{g}$ is a direct summand, $t+1\leq k\leq r$.  
The fibres look like 
\[
F_w = \bigoplus_{k=1}^t V^{(k)}\cap {}^{w}\Lie(\U )\oplus \bigoplus_{k= t+1}^r V^{(k)} \cap {}^{w}\mathfrak{u}^{(k)}
\]  
where $\mathfrak{u}^{(k)}$ is the Lie subalgebra spanned by the weights $>0$ in $\mathfrak{g}^{(k)}$.  
\[
\begin{aligned}
F_{x,xw} =\bigoplus_{k=1}^t V_k\cap {}^{x}[(\Lie(\U))\cap {}^{w}(\Lie(\U))] \oplus \bigoplus_{k= t+1}^r V_k \cap {}^{x}[\mathfrak{u}^{(k)}\cap {}^{w}\mathfrak{u}^{(k)}]
\end{aligned}
\]

\begin{lemma} Assume $J=\emptyset, \mcU=\Lie(\U)^{\oplus r}$. 
Let $x\in \W, s\in \S$ we set 
\[ h_{\overline{x}} (s):= \# \{ k\in \{1,\ldots r\} \mid x(\al_s)\in \Phi_{V^{(k)}}\}\]
where $V=\bigoplus_{k=1}^r V^{(k)}$ and $\Phi_{V^{(k)}}\subset \underline{\Phi}$ are the $T$-weights of $V^{(k)}$. 
If $x=x^i x_i$ with $x^i\in W$, then $h_{\overline{x}}(s) = h_{\overline{x_i}}(s) =: h_i(s)$. 
It holds 
\[F_{x_i} /F_{x_i,x_is} = (\mathcal{G}_{x_i(\al_s)})^{\oplus h_i(s)}.\]   
\begin{itemize}
\item[(1)] If ${}^x s\notin  W$ then 
\[h_i(s)=\# \{ k \mid V^{(k)} \subset \mcR, \; x_i(\al_s)\in \Phi_{V^{(k)}}\}. \]
\item[(2)] If ${}^xs\in  W$, then 
 \[h_i(s)=\# \{ k \mid V^{(k)} \subset \mathfrak{g},\;  x_i(\al_s)\in \Phi_{V^{(k)}}\}. \]
\end{itemize}
\end{lemma}

\paragraph{proof:}
Without loss of generality $V\subset \mcG, \mcU=\Lie(\U)$, set $x:=x_i$, we have a short exact sequence 
\[ 0\to V\cap {}^x[\mcU \cap {}^s\mcU] \to V\cap {}^x \mcU \to V\cap \mcG_{x(\al_s)}\to 0\]
Now, $V\cap \mcG_{x(\al_s)}=0$ if and only if $x(\al_s)\notin \Phi_{V}$. \\
\begin{itemize}
\item[(1)] If ${}^x s\notin  W$  then $x(\al_s)\notin \Phi$ where $\Phi$ are the $T$-weights of $\mathfrak{g}$. 
That means, if $V\subset \mathfrak{g}$ we get $h_i(s)=0$.  
\item[(2)] If ${}^xs\in W$, then $x(\al_s)\in \Phi$. This means, if $V\subset \mcR$ we get $h_i(s)=0$.
\end{itemize}
 
\hfill $\Box$

\subsection{Relative position stratification}
\subsubsection{In the flag varieties}
Let $J\subset \S$, $w\in {}^J\mathbb{W}^J$, $i, j\in I_J$. We define 
\[
\begin{aligned}
C^w &:= \mathbb{G}\phi_{e,w} \cap \left(\bigsqcup_{i\in I_J} G_i/P_i \times\bigsqcup_{i\in I_J} G_i/P_i \right)\\
C^{\leq w} &:= \overline{\mathbb{G}\phi_{e,w}} \cap \left(\bigsqcup_{i\in I_J} G_i/P_i\times \bigsqcup_{i\in I_J} G_i/P_i \right)\\
C^w_{i,j} &:= C^w\cap (G/P_i \times G/P_j)\\
C^{\leq w}_{i,j} &:= C^{\leq w}\cap (G/P_i \times G/P_j)\\
\end{aligned}
\]
For an arbitrary $w\in \W$ there exists a unique $v\in {}^J\mathbb{W}^J$ such that $\W_J w\W_J = \W_J v\W_J$, we set $ C^w:= C^v, C^w_{i,j}:= C^v_{i,j},
C^{\leq w}:= C^{\leq v}, C^{\leq w}_{i,j}:= C^{\leq v}_{i,j}$. We remark that $C^{\leq w}, C^{\leq w}_{i,j}$ are closed (but not necessary the closure of $C^w, C_{i,j}^w$, because it can happen that $C^w_{i,j}=\emptyset, C^{\leq w}_{i,j}\neq \emptyset$, see next lemma (3)).\\
Let $i,j\in I_J, \mcC_{i,j}:= \{ C_{i,j}^w \mid w\in {}^J\W^J, C_{i,j}^w\neq \emptyset \}, Orb_{i,j}:=\{G\text{-orbits in }G/P_i\times G/P_j\}$,
we have the following commutative diagram
\[
\xymatrix{
Orb_{i,j} \ar[d]_{\Phi}\ar[rrrr]^{rp} & && & W\cap {}^{x_i}\W_J \setminus  W / W\cap {}^{x_j}\W_J  \ar[d]^{\Psi} \\
\mcC_{i,j}\ar[rrrr]^{rp_{\W}} & && & \{({}^{x_i}\W_J)w({}^{x_j}\W_J)\mid w\in W\} 
}
\]
defined as follows 
\[
\begin{aligned}
 rp(G\phi_{x_i, wx_j}) &:= (W\cap {}^{x_i}\W_J) w (W\cap {}^{x_j}\W_J), \\
  rp_{\W}(C_{i,j}^w)    &:= {}^{x_i}\W_J (x_iwx_j^{-1}){}^{x_j}\W_J \\
 \Phi (G\phi_{x_i, wx_j}) &:= C_{i,j}^{x_i^{-1}wx_j}   \quad (\supset G\phi_{x_i, wx_j}) \\
 \Psi ((W\cap {}^{x_i}\W_J) w (W\cap {}^{x_j}\W_J)) &:= ({}^{x_i}\W_J )w({}^{x_j}\W_J),\\
\end{aligned}
\]
$rp, rp_{\W}$ are bijections and $\Phi, \Psi$ are surjections. 
We will from now on assume that $\Phi$, $\Psi$ are bijections as well, i.e. for every nonempty $C_{i,j}^w$ there is a $w_0\in W$ such that $\W_Jx_i^{-1} w_0 x_j\W_J = \W_J w\W_J$ and $C_{i,j}^w =G\phi_{x_i, w_0x_j} \subset G/P_i\times G/P_j$, this implies  
\[ C_{i,j}^w \cong G/( P_i \cap {}^{w_0} P_j\cap G). \]

\begin{lemma} \label{RelposInFlags} Let $J\subset \S, s\in \S\setminus{J}, i,j\in I_J$.
\begin{itemize}
\item[(1)] ${C^{\leq s}}$ is smooth, it equals $C^s\cup C^e$. 
\item[(2)] ${C^{\leq s}_{ij}} =\emptyset$ unless $Wx_j\mathbb{W}_J \in \{ Wx_i\mathbb{W}_J, Wx_is\mathbb{W}_J\}$. 
\item[(3)] Assume that $Wx_i\mathbb{W}_J \neq Wx_is\mathbb{W}_J$ and let $j \in I_J$ such that $x_isx_j^{-1} \in W$, then it holds 
\[ \iota_i (G/P_i) \neq \iota_j (G/P_j),  \quad {C^{\leq s}_{i,j}} =C^s_{i,j}, \quad {C^{\leq s}_{i,i}} =C^e_{i,i} \]
and 
$G\cap {}^{x_i}[\P_J\cap {}^{s}\P_J] = G\cap {}^{x_i} \P_{J\cap {}^s J}$, $C_{ij}^s = G/ (G\cap {}^{x_i} \P_{J\cap {}^s J}) $. 
\item[(4)]  Assume that $Wx_i\mathbb{W}_J = Wx_is\mathbb{W}_J=Wx_j\mathbb{W}_J$, then it holds $i=j$, in particular 
\[ \iota_i (G/P_i) = \iota_j (G/P_j), \quad  C^w_{i,j} =C^w_{i,i}, \text{ for all }w \]
and the first equality implies $({}^{x_i}\mathbb{P}_J )\cap G\neq ({}^{x_is}\mathbb{P}_J)\cap G$,
there is an isomorphism of $G$-varieties 
\[ G\times^{P_i} \left(({}^{x_i}\P_{J\cup \{s\}}\cap G)/P_i \right)\to {C^{\leq s}_{i,i}}, \quad (g,hP_i)\mapsto ( gP_i , ghP_i) .\]
\end{itemize}
\end{lemma}

\paragraph{proof:}
\begin{itemize}
\item[(1)] The variety $\bigsqcup_{i\in I_J} G_i/P_i$ is a smooth subvariety of $\mathbb{G}/\mathbb{P}_J$ because each $G/P_i, i\in I_J$ is smooth. It is known that $\overline{\G\phi_{e,s}}= \G \phi_{e,s} \cup \G\phi_{e,e}$ is smooth in $\mathbb{G}/\mathbb{P}_J$, therefore its intersection (i.e. pullback) is smooth in $(\mathbb{G}/\mathbb{P}_J)^H$. 
\item[(2)]  Now, ${C^{\leq s}_{i,j}} ={C^s_{i,j}} \cup {C^e_{i,j}}$ and 
${C^s_{i,j}}\neq \emptyset $ iff it contains a $T$-fixed point $\phi_{x_i, vx_j}$ for a $v\in W$, that implies $x_i^{-1}vx_j\P_J=s\P_J$, i.e. there is an $f\in \P_J$ such that $vx_jf =x_is$, therefore $f\in \P_J \cap \mathbb{W} =\W_J$ and $Wx_j \W_J =Wx_is\W_J$. Similar $ C^e_{i,j}\neq \emptyset$ iff  $Wx_j \W_J =Wx_i\W_J$. 
\item[(3)] The intersection $(G/P_i) \cap (G/P_j)$ is a $G$-equivariant subset of $\G/\P_J$, therefore it is nonempty iff it contains all $T$-fixed points $vx_i \P_J =wx_j\P_J$ with $v,w\in W$. But this is equivalent to $Wx_i\W_J =Wx_j\W_J$. \\
As we have seen before $Wx_is\W_J= Wx_j \W_J$ implies $C^e_{i,j}=\emptyset, C^s_{i,i}=\emptyset$ and therefore ${C^{\leq s}_{i,j}}=C^s_{i,j}, {C^{\leq s}_{i,i}}=C^e_{i,i}$. \\
Let 
$Wx_i\W_J\neq Wx_is \W_J$, we need to show $G\cap {}^{x_i}[\P_J\cap {}^{s}\P_J] = G\cap {}^{x_i} \P_{J\cap {}^s J}$ 
Let $\underline{\Phi}=\underline{\Phi}_+\cup \underline{\Phi}^-$ be the set of roots for $(\mathbb{G}, \B , T)$ decomposing as positive and negative roots, let $\underline{\De}_J\subset \underline{\Phi}_+$ be the simple roots corresponding to $J\subset \S$ and let $\Phi$ be the roots for $(G,T)$.
It is enough to prove that the $T$-weights on $\Lie(G\cap {}^{x_i}[\P_J\cap {}^{s}\P_J])$ equal the $T$-weights on $\Lie(G\cap {}^{x_i} \P_{J\cap {}^s J})$. 

Now,  $Wx_i\W_J\neq Wx_is \W_J$ implies ${}^{x_i}s\notin W$ or equivalently $x_i(\al_s) \notin \Phi$ where $\al_s\in \underline{\Phi}_+$ is the simple root negated by $s$.
We have the $T$-weights of $\Lie({}^{x_i}\P_J)$ are  $\{ x_i(\al)\mid \al \in \underline{\Phi}_+ \cup -\underline{\De}_J \}$, \\
the $T$-weights of  $\Lie({}^{x_is}\P_J)$ are $\{ x_i(\al)\mid \al \in \underline{\Phi}_+\setminus \{\al_s\} \cup -s(\underline{\De}_J) \cup \{ -\al_s\} \}$.\\
It follows that the $T$-weights of  $\Lie(G\cap {}^{x_i}\P_J \cap {}^{x_is} \P_J)$ are  \[
\begin{aligned} 
\{x_i(\al)&\mid \al \in (\underline{\Phi}_+ \cup [-\underline{\De}_J\cap -s(\underline{\De}_J)])\cap \Phi\} \\
= \{x_i(\al)&\mid \al \in (\underline{\Phi}_+ \cup -\underline{\De}_{J\cap s^{}J})\cap \Phi\} 
\end{aligned}
\]
and these are the $T$-weights of $\Lie( G\cap {}^{x_i} \P_{J\cap {}^s J})$.

\item[(4)] The first part is by definition. Assume $Wx_i\W_J= Wx_is \W_J$ implies $x_isx_i^{-1}= ab$ with $a\in W, b\in {}^{x_i}\W_J$. Now ${}^{x_i}\P_J \cap G$ is a parabolic subgroup of $G$ conjugated to $P_{J\cap S}$, therefore 
\[ {}^{x_is}\P_J\cap G  = ({}^{x_isx_i^{-1}}({}^{x_i}\P_J))\cap G = {}^a ({}^{x_i}\P_J)\cap G = {}^a (({}^{x_i}\P_J)\cap G)\]
and assume that this is equal ${}^{x_i}\P_J\cap G$ that implies $a \in x_i\langle J\cap S\rangle x_i^{-1}$, then $x_isx_i^{-1}=ab \in {}^{x_i}\W_J$ that implies 
$s\in J$ contradicting our assumption $s\notin J$. \\
Finally, consider the closed embedding $G\times^{P_i}\left(({}^{x_i}\P_{J\cup \{s\}}\cap G)/P_i \right) \to G\times^{P_i} G/P_i$ and compose it with the $G$-equivariant isomorphism $G\times^{P_i} G/P_i \to G/P_i \times G/P_i , \overline{(g,hP_i)}\mapsto (gP_i, ghP_i)$. The image is precisely 
$C_{i,i}^s\cup C_{i,i}^e$. 
\end{itemize}
\hfill $\Box $

\subsubsection{In the Steinberg variety}

Let $w\in \mathbb{W}^J$, $i, j\in I_J$, recall that we have a map $m_J\colon Z_J \to \mathbb{G}/\mathbb{P}_J$.    
\[
\begin{aligned}
Z^w_{i,j} &:= m_{i,j}^{-1} (C_{i,j}^w)\\ 
Z^w=Z^w_J &:= \bigsqcup_{i, j\in I_J} Z_{i,j}^w \\
Z^{\leq w}=Z^{\leq w}_J &= \bigcup_{v\leq w, v\in \mathbb{W}^J} Z_J^v\\
Z^{\leq w }_{i,j} &:= \bigcup_{v\leq w, v\in \mathbb{W}^J} Z_{i,j}^v
\end{aligned}
\]

\begin{lemma}
\begin{itemize}
\item[(a)] If $C_{i,j}^w\neq \emptyset$, the restriction 
$m_{i,j}\colon Z^w_{i,j}\to C^w_{i,j}$ is a vector bundle with fibres isomorphic to $F_i \cap {}^{x_iwx_j^{-1}}F_j$, it induces a bijection on $T$-fixed points. In particular, all nonempty $Z^w_{i,j}$ are smooth.  
\item[(b)] For any $s\in \S$ the restriction $m\colon \overline{Z^s}\to {C^{\leq s}}$ is a vector bundle over its image, 
in particular $\overline{Z^s}$ is smooth. More precisely, it is a disjoint union $\overline{Z_{i,j}^s} \to {C_{i,j}^{\leq s}}$ with  
\begin{itemize}
\item[(1)]
$\overline{Z_{i,j}^s}\neq \emptyset $ implies $Wx_j\mathbb{W}_J = Wx_is\mathbb{W}_J$.
\item[(2)]   Assume that $Wx_i\mathbb{W}_J \neq Wx_is\mathbb{W}_J$, then 
$\overline{Z_{i,j}^s}=Z_{i,j}^s$ and $\overline{Z_{i,i}^s}=\emptyset $.
\item[(3)] Assume that $Wx_i\mathbb{W}_J = Wx_is\mathbb{W}_J$, then it holds  
$\overline{Z_{i,i}^s} \to {C_{i,i}^{\leq s}}$ is a vector bundle. 
\end{itemize}
\end{itemize}
\end{lemma}

\paragraph{proof:}
\begin{itemize}
\item[(a)] 
As $C^w_{i,j}$ is assumed to be a diagonal $G$-orbit in $G/P_i\times G/P_j$, it is a homogeneous space and the statement easily follows from a wellknown lemma, 
cp. \cite{Sl}, p.26, lemma 4.
\item[(b)] 

\begin{itemize}
\item[(1)] If $\overline{Z_{i,j}^s}\neq \emptyset$, then $C_{i,j}^s\neq \emptyset$ and by the proof of the previous lemma \ref{RelposInFlags}, (2), the claim follows.
\item[(2)] If $Wx_i\mathbb{W}_J \neq Wx_is\mathbb{W}_J$, then by lemma \ref{RelposInFlags}, (3), $C_{i,j}^{\leq s}=C_{i,j}^s$ is already closed, therefore 
$Z_{i,j}^s$ is closed as well. Also, $C_{i,i}^{\leq s}=C_{i,i}^e$ is already closed, therefore $Z_{i,i}^e$ is closed as well.  
\item[(3)] If $Wx_i\mathbb{W}_J = Wx_is\mathbb{W}_J$, then $C_{i,i}^{\leq s} $ is the closure of the $G$-orbit $C_{i,i}^s$ and by lemma \ref{RelposInFlags}, (4) we have $G\times^{P_i} \left(({}^{x_i}\P_{J\cup \{s\}}\cap G)/P_i \right)\to {C^{\leq s}_{i,i}}, \quad (g,hP_i)\mapsto ( gP_i , ghP_i)$ is an isomorphism. 
We set $X:=$
\[
\{ (gf, gP_i,ghP_i) \in G(F_i \cap {}^{{}^{x_i}s}F_i) \times G/P_i\times G/P_i 
\mid g\in G, f\in F_i \cap {}^{x_isx_i^{-1}}F_i, h \in {}^{x_i}\P_{J\cup \{s\}}\cap G \} 
\] 
and we claim $\overline{Z_{i,i}^s}=X$. First, observe that $X\subset Z_{i,i}$ because $gf=gh (h^{-1}f)$ with $h^{-1}f \in F_i \cap {}^{x_isx_i^{-1}}F_i$. 
One can easily check the following steps.
\begin{itemize}
\item[(*)] $X\to C_{i,i}^{\leq s}$ is a vector bundle with fibre over $F_i \cap {}^{x_isx_i^{-1}}F_i$. In particular, we get that $X$ is smooth irreducible and $\dim X=\dim Z_{i,i}^s$.
\item[(*)] $Z_{i,i}^s\subset X$. 
\item[(*)] $X$ is closed in $Z_{i,i}$ because we can write it as $X=p^{-1}(G (F_i \cap {}^{x_isx_i^{-1}}F_i)) \cap m^{-1}(C_{i,i}^{\leq s})$. Since 
$F_i \cap {}^{x_isx_i^{-1}}F_i$ is (by definition) $B_i={}^{x_i}B$-stable, we get $G(F_i \cap {}^{x_isx_i^{-1}}F_i)$ is closed in $V$. This implies $X$ is closed. 
\end{itemize}
\end{itemize}
\end{itemize}

\section{A short lamentation on the parabolic case}
From the next section on we assume that all $P_i=B_i$ are Borel subgroups. What goes wrong with the more general assumption (which we call the \emph{parabolic case})?
\begin{itemize}
\item[(1)] We do not know whether $C^w_{i,j}$ (see previous section) is always a $G$-orbit. That is relevant for Euler class computation in Lemma 14.  
\item[(2)] The cellular fibration property has to be generalized because 
$C^w:=\{gP, gwP^\prime \mid g\in G\}\subset G/P \times G/P^{\prime} \xrightarrow{ pr_1} G/P$ is not a vector bundle (its fibres are unions of Schubert cells). This complicates Lemma 12.
\item[(3)] We do not know what is the analogue of lemma 13, i.e. what can we say about 
$\mcZ^{\leq x}*\mcZ^{\leq y}$ ?
\item[(4)] The cycles $[\overline{Z_{i,j}^s}]$ are not in general multiplicative generators. If we try to understand more generally $[\overline{Z_{i,j}^w}]$, 
the multiplicity formular does not give us as much information as for $[\overline{Z_{i,j}^s}]$ because $\overline{Z_{i,j}^s}$ is even smooth. 
Also understanding the $[Z_{i,j}^w]$ is not enough, since they do not give a basis as a free $\mcE$-module because the rank is wrong (cp. failing of cellular fibration lemma). 
\end{itemize}

The point (4) is the biggest problem. Even for $H_*^G(G/P\times G/P)$ we do not know a set of generators and relations (see next chapter). \\
\textbf{So, from now on we assume $J=\emptyset$}.

\subsection{Convolution operation on the equivariant Borel-Moore homology of the Steinberg variety}

\begin{defi} Let $H\in \{ pt, T,  G\}$ with $T\subset  G$ where $T$ is a maximal torus. \\
We define the \emph{$H$-equivariant algebra of a point} to be $H_H^*(pt)$ with product equals the cup-product, we will always identify it with $H^H_*(pt):= H_H^{-*}(pt)$. It is a graded $\C$-algebra concentrated in negative even degrees.\\
We define the \emph{$H$-equivariant Steinberg algebra} to be the $H$-equivariant  Borel-Moore homology algebra of the  Steinberg variety, the product is the convolution product, see \cite{CG}, \cite{VV}.\\
We say ($H$-equivariant)\textbf{company algbra} to the $H$-equivariant cohomology algebra of $E$, the product is the cup-product.  
\[ 
\begin{aligned}
\La_H &:= H_H^*(pt) \text{ for the }H\text{-equivariant algebra of a point}, \\
\mcZ_H &:= H^H_*(Z) \text{ for the }H\text{-equivariant Steinberg algebra}, \\
\mcE_H &:= H_H^* (E) \text{ for the }H\text{-equivariant company algebra} . 
\end{aligned}
\]
For $H=pt$ we leave out the adjective $H$-equivariant and leave out the index $H$. 
\end{defi}
Recall, that $\mcZ_H$ and $\mcE_H$ are left graded modules over $\La_H$. 
Furthermore, $\mcE_H$ is a left module over $\mcZ_H$. This follows from considering $M_1=M_2=M_3=E$ smooth manifolds and 
$Z\subset M_1\times M_2, E=E\times \{ \overline{(e,0)}\} \subset M_2\times M_3$. Then the 
set-theoretic convolution gives $Z\circ E=E$, which implies the operation. \\

Also, $\mcZ_H$ is a left module over $\mcE_H$. This follows from considering $M_1=M_2=M_3=E$ smooth manifolds ($\dim_{\C}E=:e$) and 
$E\inj M_1\times M_2$ diagonally, $Z\subset M_2\times M_3$, then the set-theoretic convoltion gives $E\circ Z=Z$, that implies that we have a map
\[ 
H_{2e_i-p}^H(E_i)\times H_{e_i+e_j-q}^H(Z_{i,j}) \to H_{e_i+e_j-(p+q)}^H(Z_{i,j})
\]
Using Poincare duality we get $H_{2e_i-p}^H(E)\cong H_H^{p}(E)$ and the grading $H_{[q]}^H(Z):= \bigoplus_{i,j}H_{e_i+e_j-q}^H(Z)$ the previous map gives an 
operation of the $H_H^*(E)$ on $H_{[*]}^H(Z)$ which is $H^H_*(pt)$-linear. 
We denote the operations by 
\[ 
\begin{aligned}
*\colon &\mcZ_H \times \mcE_H &\to \mcE_H \\
\diamond \colon &\mcE_H \times \mcZ_H &\to \mcZ_H
\end{aligned}
\] 
Furthermore, there are forgetful algebra homomorphisms 
\[
\begin{aligned}
\La_G &\to \La_T &\to  \La, \\ 
\mcZ_G &\to \mcZ_T &\to \mcZ, \\ 
\mcE_G &\to \mcE_T &\to \mcE. 
\end{aligned}
\] 

Let us investigate some elementary properties of the convolution operations. 
From \cite{VV2}, section 5, p.606, we know that the operation of $\mcZ_{G}$ on $\mcE_G$ is faithful, i.e. we get an injective $\C$-algebra homomorphism 
 \[ \mcZ_G \inj \End_{\C -alg} (\mcE_G ).\]
We have the following cellular fibration property. 
We choose a total order $\leq $ refining Bruhat order on $\W$. For each $i,j\in I$ we get a filtration into closed $G$-stable subsets of $Z_{i,j}$ 
by setting $Z_{i,j}^{\leq w}:=\bigcup_{v\leq w} Z_{i,j}^v, \; w\in \W$. Via the first projection $pr_1\colon C_{i,j}^v\to G/B_i$ is a $G$-equivariant vector bundle with fibre $B_ivB_j/B_j$, we call its (complex) dimension $d_{i,j}^v$, also $Z_{i,j}^v\to C_{i,j}^v$ is a $G$-equivariant vector bundle, we define the complex fibre dimension $f_{i,j}^v$. By the $G$-equivariant Thom isomorphism (applied twice) we get 
\[ H_m^G(Z_{i,j}^v)= H_{m-2d_{i,j}^v-2f_{i,j}^v}^G(G/B_i). \]
In particular, it is zero when $m$ is odd and $H_*^G(Z_{i,j}^v)$  is a free $H_*^G(pt)$-module with basis \\
$b_x, \; x\in W, \deg b_x= 2 \dim (B_ixB_i)/B_i + 2d_{i,j}^v+2f_{i,j}^v$. \\ 
Using the long exact localization sequence in $G$-equivariant Borel-Moore homology for every $v\in \W$, we see that 
$Z_{i,j}^{v}$ is open in $Z_{i,j}^{\leq v}$ with an closed complement $Z_{i,j}^{<v}$. 
We conclude inductively using the Thom isomorphism that $H_{odd}^G(Z_{i,j}^{\leq v})=0$ and that $H_*^G(Z_{i,j}^{\leq w})=\bigoplus_{v\leq w} H_*^G(Z_{i,j}^v)$. We observe, that 
$\# \{ w\in \W\mid Z_{i,j}^w\neq \emptyset \} =\# W$ for every $i,j\in I$. It follows that $H_*^G(Z_{i,j})$ is a free 
$H_*^G(pt)$-module of rank $\# (W\times W)$, and that every $H_*^G(Z_{i,j}^{\leq v})\xrightarrow{i_*} H_*^G(Z_{i,j})$ is injective. \\
We can strengthen this result to the following lemma.

\begin{lemma} Let $\leq$ be a total order refining Bruhat order on $\W$. 
For any $w\in \W$ set $Z^{\leq w} :=m^{-1} ( \bigcup_{v\leq w}C^v  ) =\bigcup_{v\leq w } Z^v$. 
The closed embedding $i\colon Z^{\leq v} \to Z$ gives rise 
to an injective morphism of $H^*_G(E)$-modules $i_*\colon \mcZ^{\leq v}_G := H^G_*(Z^{\leq v})\to \mcZ_G $. We identify in the following $\mcZ^{\leq v}_G$ with its image in $\mcZ_G$. For all $v\in W$ we have 

\[ 
\begin{aligned}
\mcZ_G^{\leq w}  &=\bigoplus_{v\leq w} \mcE_G \diamond [ \overline{Z^v}]  \quad & \text{ as } \mcE_G\text{-module}\\
1_i *\mcZ_G^{\leq w} * 1_j &=\bigoplus_{v\leq w} \mcE_i \diamond [ \overline{Z^v_{i,j}}] \quad & \text{ as } \mcE_i\text{-module}
\end{aligned}
\] 
where $\mcE_i =H_G^*(E_i)$. Each $[\overline{Z^v}]$ is nonzero (and not necessarily a homogeneous element). 
In particular, $\mcZ_G$ (as ungraded module) is a free left $\mcE_G$-module of rank $\# \W$.
\end{lemma}
 
\paragraph{proof:}
Now first observe that set-theoretically we have $E\circ Z^v=Z^v$ (where we use the diagonal embedding for $E$ again). This implies that the direct sum decomposition $H_*^G(Z)=\bigoplus_{v\in \W} H_*^G(Z^v)$ is already a decomposition of $H_G^*(E)$-modules. 

Now we know that we have by the Thom-isomorphism algebra isomorphims  
\[H^*_G(E)\cong H^*_G(\bigsqcup_{i\in I} G/B_i) \cong H^*_G(Z^v),\]
using that $\#\{(i,j) \mid  Z_{i,j}^v\neq \emptyset \}= \#I$. 
Now, Poincare duality is given by $H_G^q(Z^v_{i,j}) \to H_{2\dim Z^v_{i,j}-q}^G (Z^v_{i,j})$, $ \al \mapsto \al \cdot [Z^v_{i,j}]$ the composition gives 
\[ H_G^p(E_i) \to H_{2\dim Z^v_{i,iv}-q}^G (Z^v_{i,iv}),\quad  c \mapsto  c \cdot [Z^v_{i,iv}] .\]
\hfill $\Box $ 
 
\begin{lemma} \label{ConvolutionOfSubmodules}
For each $x,y\in \W$ with $l(x)+l(y)=l(xy)$ we have 
\[ \mcZ^{\leq x}_G *\mcZ^{\leq y}_G \subset \mcZ^{\leq xy}_G \]  
\end{lemma}

\paragraph{proof:}
By definition of the convolution product, it is enough to check that for all $w\leq x, v\leq y$ it holds for the set theoretic convolution product 
\[ {Z^w_{i,j}}\circ {Z^v_{j^\prime,k}} \subset \begin{cases} \emptyset , &\quad  j\neq j^\prime \\ 
 Z^{\leq xy}_{i,k}, &\quad  j=j^\prime 
 \end{cases}
  \]
for $i,j,j^\prime,k\in I$, 
because by definition $Z^{\leq x} \circ Z^{\leq y} = \bigcup_{w\leq x, v\leq y} Z^w\circ Z^v$. Now, the case $j\neq j^\prime$ follows directly from the 
definition. Let $j=j^\prime$.  
Let $\C^w:= \G(\B, w\B) \subset \G/\B\times \G/\B$. According to Hinrich, Joseph \cite{HJ}, 4.3 it 
holds $\C^w \circ \C^v \subset {\C^{wv}}$ for all $v,w\in \W$. Now, we can adapt this argument to prove that $C_{i,j}^w\circ C_{j,k}^v \subset 
C_{j,k}^{wv}$ as follows: \\
Since $C_{i,j}^w\neq \emptyset, C_{j,k}^v \neq \emptyset$ we have that $w_0=x_iwx_j^{-1}\in W, v_0=x_j vx_k^{-1}\in W$ and \\
$C_{i,j}^w=G(B_i, w_0B_j), C_{j,k}^v=G(B_j,v_0B_k)$. We pick $M_1=G/B_i, M_2=G/B_j, M_3=G/B_k$ for the convolution and get 
\[p_{13}( p_{12}^{-1}C_{i,j}^w\cap p_{23}^{-1}C_{j,k}^v) = \{ g(B_i, w_0bv_0B_k)\mid g\in G, b\in B_j\}.\] 
Now since the length are adding one finds  
$B_iw_0 B_j v_0B_k = B_i (w_0v_0)B_k$ , as follows 
\[
\begin{aligned}
   w_0 B_j v_0B_k & = x_i[w ({}^{x_j^{-1}}G\cap \B) v ({}^{x_k^{-1}}G\cap \B)]x_k^{-1} \\
  & \subset x_i[ w\B v\B ]x_k\cap G \subset x_i[\B wv \B]x_k^{-1} \cap G   \\
  &= [{}^{x_i}\B (x_iwvx_k^{-1}) {}^{x_k}\B] \cap G
  &= B_i w_0v_0 B_k
  \end{aligned}
  \]
For the last equality, clearly $  B_i w_0v_0 B_k \subset  [{}^{x_i}\B (x_iwvx_k^{-1}) {}^{x_k}\B] \cap G$. Assume $[{}^{x_i}\B (x_iwvx_k^{-1}) {}^{x_k}\B] \cap G= \bigcup B_itB_k$ for certain $t\in W$, then clearly $B_itB_k\subset [{}^{x_i}\B (x_iwvx_k^{-1}) {}^{x_k}\B] \cap G \cap [{}^{x_i}\B t {}^{x_k}\B] \cap G$ as this intersection is empty if $t\neq (x_iwvx_k^{-1})$, the last equality follows. \\
Then using $Z_{i,j}^w=\{ g(f_i=w_0f_j, B_i, w_0B_j) \in V \times G/B_i\times G/B_j \mid g\in G, f_i\in F_i, f_j\in F_j\}$ one concludes by definition that 
$Z_{i,j}^w\circ Z_{j,k}^v \subset Z_{j,k}^{wv}$
\hfill $\Box$ \\

We have the following corollary whose proof we have to delay until we have introduced the localization to the $T$-fixed point.
\begin{coro} \label{ConvOfZs}
For $s\in \S, w\in \W$ with $l(sw)= l(w)+1$, 
\[ [\overline{Z^s}] * [\overline{Z^w}] =[\overline{Z^{sw}}] \text{  in  } \mcZ_G^{\leq sw}/\mcZ_G^{< sw}.\] 
Since $[\overline{Z^v}]= \sum_{s,t\in I} [\overline{Z_{s,t}^v}]$ for all $v\in \W$, this is equivalent to 
$i,j,l,k\in I$ we have 
\[ [\overline{Z^s_{i,j}}] * [\overline{Z^w_{l,k}}] =\delta_{l,j}[\overline{Z^{sw}_{i,k}}] \text{  in  } \mcZ_G^{\leq sw}/\mcZ_G^{< sw}.\] 
\end{coro}


\subsubsection{Computation of some Euler classes}

\begin{defi} (Euler class)
 Let $M$ be a finite dimensional complex $\mathfrak{t}=Lie(T)$-represenation. Then, we have a weight space decomposition
\[ M =\bigoplus_{\al \in \Hom_{\C}(\mathfrak{t}, \C) } M_{\al}, \quad M_{\al} = \{ m\in M \mid tm =\al(t) m \}. \] 
We define 
\[ 
\eu(M):= \prod_{\al\in \Hom (\mathfrak{t}, \C)} \al^{\dim M_{\al}} \quad \in \C[\mathfrak{t}] =H_T^*(pt) 
\]
For a $T$-variety $X$ and a $T$-fixed point $x\in X$, we define the \textbf{Euler class} of $x\in X$ to be
\[ 
\eu (X,x):= \eu( T_x X),
\] 
where the $\mathfrak{t}$-operation on the tangent space $T_xX$ is the differential of the natural $T$-action. \\
Observe, that $\eu (T_x^*X)=(-1)^{\dim T_xX} \eu (T_xX)$.
\end{defi}
%
%

Recall from an earlier section the notation $Z^w := m^{-1} (C^w)$.  
We are particularly interested in the following Euler classes, let $w=w^k x_k, x=x^i x_i, y= y^j x_j\in \W, w^k, x^i, y^j\in W$  
\[
\begin{aligned}
\La_w &:= \eu ( E, \phi_w ) = \eu ( T_{\phi_{w^kx_k}} E_k ) ,   \quad & \in H_T^*(pt) \\
eu (\overline{Z^w}, \phi_{x,y}) &= (\eu ( T_{\phi_{x^ix_i, y^jx_j}} \overline{Z_{ij}^w}))^{-1} , \quad  & \in K:=Quot(H_T^*(pt))
\end{aligned}
\]

Remember $F_w:=\mu^{-1}(\phi_{w})= \mu_k^{-1}(\phi_{w^{k}x_k}) = {}^{w^k}F_k, \;\; 
F_{x,y}:= m^{-1} (\phi_{x,y}) = {}^{x^i}F_i \cap {}^{y^j}F_j = F_x\cap F_y$. In particular, we can see them as $\mathfrak{t}$-representations. We also consider the following $\mathfrak{t}$-representations 
\[
\begin{aligned}
\mathfrak{n}_w& := T_{w^kP_k} G/P_k = \mathfrak{g}\cap {}^{w}\mcU^-={}^{w^k}[ \mathfrak{g}\cap {}^{x_k}\mcU^-] \\
\mathfrak{m}_{x,y}& :=  \frac{\mathfrak{n}_x}{ \mathfrak{n}_x\cap \mathfrak{n}_y} = 
\mathfrak{g}\cap \frac{{}^{x}\mcU^-}{{}^{x}\mcU^-\cap {}^{y}\mcU^-}  
\end{aligned}
\] 
where $\mcU^-:=\Lie(\U^-)$ with $\U^-\subset \B^-:={}^{w_0}\B$ is the unipotent radical where $w_0\in \W$ is the longest element. 
Some properties can easily be seen. 
\begin{itemize}
\item[(1)] $\mathfrak{n}_x =\prod_{\al \in \Phi \cap x^{-1}\underline{\Phi}^-} \al$. 
\item[(2)] If $s\in \S$, $x\in \W$ such that ${}^x s\in W$, then 
\[ \eu (\mathfrak{n}_x ) = -\eu (\mathfrak{n}_{xs} ), \quad \eu (\mathfrak{m}_{x,xs}) = -\eu(\mathfrak{m}_{xs,x})= x(\al_s) \]
\item[(3)] If $s\in \S$, $x\in \W$ such that ${}^xs\notin W$, then 
\[ \mathfrak{n}_x  =  \mathfrak{n}_{xs} ,  \quad \eu (\mathfrak{m}_{x,xs}) = \eu(\mathfrak{m}_{xs,x})= 0 \] 
\end{itemize}
Furthermore, for $s\in S, x\in \W, i\in I$ we write set as a shortage
\[
\begin{aligned}
Q_{x}(s) &:= \eu(F_x/F_{x,xs}),\\
Q_i(s) &:= Q_{x_i}(s), \\
q_i(s) &:= \prod_{\al \in \Phi_{\mcU}, s(\al) \notin \Phi_{\mcU}, x_i(\al) \in \Phi_V} \al .
\end{aligned}
\]
for $x=x^ix_i$ with $x^i\in W$ it holds $Q_x(s)=x^i(Q_i(s)), Q_i(s)=x_i(q_i(s))$, i.e. 
\[ Q_x(s) =x(q_i(s)) \]

\begin{lemma} Let $J=\emptyset$, it holds 
\begin{itemize}
\item[(1)] for $w\in \W$
\[ 
\La_w = \eu (F_w \oplus \mathfrak{n}_w) 
\]
\item[(2)] If $s\in \S , x\in \W, \al_s\in \underline{\Phi}^+$ with $s(\al_s)=-\al_s$ and ${}^x s\in W$
\[
\begin{aligned}
\eu (\overline{Z^s}, \phi_{x,xs}) &= \eu (F_{x,xs}\oplus \mathfrak{n}_x\oplus \mathfrak{m}_{x,xs})= \;x(\al_s)\;Q_x(s)^{-1}\;\La_x  \\
 \eu (\overline{Z^s}, \phi_{x,x}) &= \eu ( F_{x,xs} \oplus \mathfrak{n}_x\oplus \mathfrak{m}_{xs,x})= \;-\eu (\overline{Z^s}, \phi_{x,xs}). 
 \end{aligned}
 \]
\item[(3)] If $s\in \S , x\in \mathbb{W}$ and ${}^xs\notin W$
 \[ \eu (\overline{Z^s}, \phi_{x,xs}) = \eu (F_{x,xs}\oplus \mathfrak{n}_{x}) = Q_x (s)^{-1}\La_{x}  \]
 \item[(4)] Let $x,w\in \W$. Then 
 \[ 
  \eu (\overline{Z^w}, \phi_{x,xw}) = \eu(F_{x,xw}\oplus \mathfrak{n}_x\oplus \mathfrak{m}_{x,xw})
 \]

\end{itemize}
\end{lemma}

\paragraph{proof}
\begin{itemize}
\item[(1)]
We know $\mu_k \colon E_k\to G/B_k,\; B_k=G\cap {}^{x_k}\B $ is a vector bundle, therefore we have a short exact sequence of tangent spaces 
\[ 
0 \to T_{\phi_{w}} \mu_k^{-1}(w^kB_k) \to T_{\phi_w} E_k  \to T_{w^kB_k} G/B_k \to 0 
\]
which is a split sequence of $T$-representations implying the first statement. 
\item[ad (3,2)] 
Let $i,j\in I_J$ such that $x^i:=xx_i^{-1},y^j:=xsx_j^{-1}\in W$. \\

\item[(2)] If ${}^xs\in W$ we have that $i=j$
and $\overline{Z_{i,i}^s} \to C_{i,i}^{\leq s}\cong G\times^{B_i} (G\cap {}^{x_i} \P_{\{s\}} )/B_i$ is a vector bundle. For $x^\prime \in \{ x,xs\}$ we have a short exact sequence on tangent spaces  
\[ 0 \to F_{x, xs}\to T_{\phi_{x,x^\prime}} \overline{Z_{i,i}^s} \to T_{\phi_{x,x^\prime}} C_{i,i}^{\leq s} \to 0 \]

Using the isomorphism $G\times^{B_i} [({}^{x_i}\P_{\{s\}} \cap G)/B_i] \to C_{i,i}^{\leq s}, \; (g,hB_i) \mapsto (gB_i, ghB_i)$ we get 
\[ \eu( T_{\phi_{x,x^\prime}} C_{i,i}^{\leq s}  ) = \begin{cases} \eu (T_{\overline{(x^i, B_i)}} G\times^{B_i} [({}^{x_i}\P_{\{s\}} \cap G)/B_i]) 
= \eu (\mathfrak{n}_x)\cdot \eu(\mathfrak{m}_{xs,x}), \; \; x^\prime =x \\
\eu (T_{\overline{(x^i, {}^{x_i}sB_i)}} G\times^{B_i} [({}^{x_i}\P_{\{s\}} \cap G)/B_i]) 
= \eu (\mathfrak{n}_x)\cdot \eu (\mathfrak{m}_{xs,x}), \; \; x^\prime =xs 
\end{cases}
\]
It follows $\eu (\overline{Z^s}, \phi_{x,x})= \eu (F_{x,xs}) \cdot \eu (\mathfrak{n}_x)\cdot \eu (\mathfrak{m}_{xs,x})$ and $\eu (\overline{Z^s}, \phi_{x,xs})= \eu 
(F_{x,xs}\oplus \mathfrak{n}_x\oplus \mathfrak{m}_{x,xs})$.
\item[(3)] If ${}^xs\notin W$ we get $i\neq j$ and $Z_{i,j}^s$ is closed and a vector bundle over 
$C_{i,j}^s = G/ (G\cap {}^x \B)$, we get a short exact sequence on tangent spaces 
\[ 0\to F_{x,xs} \to T_{\phi_{x,xs}} Z_{i,j}^s \to T_{\phi_{x,xs}} C_{i,j}^s \to 0.\]
We obtain $\eu (\overline{Z^s}, \phi_{x,xs})= \eu (F_{x,xs}) \eu ( \mathfrak{n}_{x})$. 
\item[(4)]
Pick $i,j\in I$ such that $x\in Wx_i, xw\in Wx_j$.
We have the short exact sequence 
\[
0\to F_{x,xw} \to T_{\phi_{x,xw}} \overline{Z_{i,j}^w}\to T_{\phi_{x,xw}}C_{i,j}^w \to 0
\] 
Then, recall the isomorphism 
\[
\begin{aligned}
C_{i,j}^w =G\phi_{x,xw} &\to G/(G\cap{}^x\B\cap {}^{xw}\B) \\
\phi_{x,xw} &\mapsto \overline{e}:=e(G\cap{}^x\B\cap {}^{xw}\B)
\end{aligned}
\] 
Again we have a short exact sequence
\[ 
0\to T_{\overline{e}} (G\cap {}^x\B)/(G\cap{}^x\B\cap {}^{xw}\B) \to T_{\overline{e}} G/(G\cap{}^x\B\cap {}^{xw}\B) \to T_{\overline{e}} G/(G\cap{}^x\B) \to 0 
\]  
Together it implies $\eu (\overline{Z_{i,j}^w}, \phi_{x,xw}) = \eu (F_{x,xw})\eu (\mathfrak{n}_x / (\mathfrak{n}_x\cap \mathfrak{n}_{xw})) \eu (\mathfrak{n}_x). $
\end{itemize}

\hfill $\Box $

\begin{coro} Let $J=\emptyset, \mcU=\Lie(\U)^{\oplus r}$, it holds 
\begin{itemize}
\item[(1)] If $s\in \S , x\in \mathbb{W}$ and ${}^x s\in W$, then $h_{\overline{x}}(s)=h_{\overline{xs}}(s)$ and 
\[
\begin{aligned}
\La_{x} &= (-1)^{1+h_{\overline{xs}}(s)}\La_{xs} \\
\eu (\overline{Z^s}, \phi_{x,xs}) &= (x(\al_s))^{1-h_{\overline{x}}(s)}\La_x 
 \end{aligned}
 \]
\item[(2)] If $s\in \S , x\in \mathbb{W}$ and ${}^xs\notin W$
 \[ \eu (\overline{Z^s}, \phi_{x,xs})= x(\al_s))^{-h_{\overline{x}}(s)}\La_x \]
\end{itemize} 

\end{coro}

\paragraph{proof:} This follows from $q_x(s)= x(\al_s)^{h_{\overline{x}}(s)}$ and \\
if ${}^xs\in W$ we have that $i=j$ and $h_{\overline{x}}(s)=h_{\overline{xs}}(s)$. Therefore we get 
\[
\begin{aligned}
\eu (F_x) & = x(\al_s)^{h_{\overline{x}}(s)} \eu (F_{x,xs}) \\
       &= (-1)^{h_{\overline{xs}}(s)} (xs(\al_s))^{h_{\overline{xs}}(s)} \eu (F_{xs,x}) \\
       &= (-1)^{h_{\overline{xs}}(s)} \eu (F_{xs})
\end{aligned}
\] 
Using that $\eu(\mathfrak{n}_x)=-\eu(\mathfrak{n}_{xs})$ we obtain $\La_x = (-1)^{1+h_{\overline{xs}}(s)}\La_{xs}$
\hfill $\Box$


\subsubsection{Localization to the torus fixed points}
Now, we come to the application of localization to $T$-fixed points.
We remind the reader that $Z$ is a cellular fibration and $E$ is smooth, therefore in both 
cases the odd ordinary (=singular) cohomology groups vanish for $Z$ and $E$. This implies in 
particular that $E,Z$ are \emph{equivariantly formal}, which is (in the case of finitely $T$-fixed points) equivalent to $\mcZ_G$ and $\mcE_G$ are free 
modules over $H_G^*(pt)$.

If we denote by $K$ the quotient field of $H_G^*(pt)$ and for any $T$-variety $X$ 
\[ H^T_*(X) \to \mcH_* (X):= H_*^T(X) \otimes_{H_T^*(pt)} K ,\quad  \al \mapsto \al \otimes 1 .\]

\begin{lemma}
\begin{itemize}
\item[(1)]
\[ \mcH_* (E) = \bigoplus_{w\in \W} K \psi_w, \quad \mcH_* (Z) =\bigoplus_{x,y\in \W} K \psi_{x,y} \]
where $\psi_w=[\{\phi_w\}]\otimes 1, \psi_{x,y}=[\{\phi_{x,y}\}]\otimes 1$. 
\item[(2)] For every $i\in I$, $w\in Wx_i$ we have a map $w\cdot \colon \mcE_i:=H_G^*(E_i)\to \C[\mathfrak{t}]$, via taking the forgetful map composed with 
the pullback map under the closed embedding $i_w\colon \{\phi_w\} \to E_i$ 
\[
 \mcE_i =H_G^*(E_i)\to H_T^*(E_i) \xrightarrow{i_w^*} H_T^*(pt)=\C[\mathfrak{t}],  
\]
we denote the map by $f\mapsto w(f)$, $f\in \mcE_i, w\in \W$. 
Furthermore, composing the forgetful map with the map from before we get an injective algebra homomorphism 
 \[
 \begin{aligned}
   \Theta_i\colon \mcE_i &\to H_T^*(E_i)\inj H_T^*(E_i)\otimes K\cong \bigoplus_{w\in Wx_i} K\psi_w\\ 
   c& \quad \longmapsto \quad \sum_{w\in Wx_i} w(c) \La_{w}^{-1} \psi_w.
   \end{aligned}
   \]
We set $\Theta=\bigoplus_{i\in I} \Theta_i\colon \mcE_G \to \bigoplus_{w\in \W} K\psi_w$.   
\end{itemize}
\end{lemma}

\paragraph{proof:}
\begin{itemize}
\item[(1)]
This is GKM-localization theorem for $T$-equivariant cohomology, for a source also mentioning the GKM-theorem for 
$T$-equivariant Borel-Moore homology see for example \cite{Br}, Lemma 1. 

\item[(2)]
This is \cite{EG2}, Thm 2, using the equivariant cycle class map to identifiy $T$-equivariant Borel-Moore homology of $E$ with the $T$-equivariant Chow ring. 
\end{itemize}
\hfill $\Box$

\paragraph{The $\W$-operation on $\mcE_G\colon$} Recall that the ring of regular functions $\C[\mathfrak{t}]$ on $\mathfrak{t}=Lie(T)$ is a left $W$-module and a left $\W$-module with respect to $ w\cdot f (t) = f(w^{-1}t w), \; w\in \W(\supset W)$. 
The from $W$ to $\W$ induced representation is given by 
\[ 
\Ind_{W}^{\W} \C[\mathfrak{t}] = \bigoplus_{i\in I} x_i^{-1} \C[\mathfrak{t}],
\]  
for $w\in \W, i\in I$ the operation of $w$ on $x_i^{-1}\C[\mathfrak{t}]$ is given by
\[
\begin{aligned}
x_i^{-1} \C[\mathfrak{t}] &\to x_{iw^{-1}}^{-1}\C[\mathfrak{t}]\\
x_i^{-1} f &\mapsto wx_i^{-1} f
\end{aligned}
\]  
where we use that $wx_i^{-1}W= x_{iw^{-1}}^{-1}W$. \\
Now, we identify $\mcE_G=\bigoplus_{i\in I} \mcE_i$ with the left $\W$-module $\Ind_{W}^{\W} \C[\mathfrak{t}]$ via $\mcE_i=x_i^{-1}\C[\mathfrak{t}]$.

Furthermore, we have the (left) $\W$-representation on $\bigoplus_{x\in \W} K(\La_x^{-1}\psi_{x})$ defined via 
\[ w(k (\La_x^{-1} \psi_x)) := k(\La_{xw^{-1}} \psi_{xw^{-1}}), \quad k\in K, w\in \W. \]

\begin{lemma} The map $\Theta\colon \mcE_G \to \bigoplus_{x\in \W} K(\La_x^{-1}\psi_{x})$ is $\W$-invariant. 
\end{lemma}

\paragraph{proof:} Let $w\in \W$, we claim that there is a commutative diagram 
\[
\xymatrix{ \mcE_G \ar[rrr]^{\Theta}\ar[d]_{w\cdot} &&& \bigoplus_{x\in \W} K(\La_x^{-1}\psi_{x}) \ar[d]_{w}\quad  & c \ar[r]\ar[d] &\sum_{x\in \W} i_x^*(c)\ar[d] \La_x^{-1} \psi_x  \\
\mcE_G \ar[rrr]^{\Theta} &&& \bigoplus_{x\in \W} K(\La_x^{-1}\psi_{x})\quad  &w\cdot c \ar[r] &\sum_{x\in \W} i_{xw}^*(c) \La_x^{-1} \psi_x
}
\]
We need to see $i_x^*(w\cdot c) = i_{xw}^*(c)$.  Let $xw\in Wx_i$, $x\in Wx_{iw^{-1}}$
This means that the diagram 
\[
\xymatrix{
\mcE_i\ar[dr]_{i_{xw}^*} \ar[rr]^{w\cdot } & &\mcE_{iw^{-1}}\ar[dl]_{i_x^*} \\ 
& H_T^*(pt) & 
}
\] 
is commutative. But it identifies with 
\[
\xymatrix{
x_i^{-1} \C[\mathfrak{t}]\ar[dr]_{{xw}\cdot } \ar[rr]^{w\cdot }& & x_{iw^{-1}}^{-1} \C[\mathfrak{t}] \ar[dl]_{x\cdot } & &  x_i^{-1} f \ar[rr]\ar[dr] &&  wx_i^{-1} f \ar[dl] \\
& \C[\mathfrak{t}] & && &  xwx_i^{-1} f. & 
}
\]
The diagram is commutative.  
\hfill $\Box $ 

\begin{rem}
From now on, we use the following description of the $\W$-operation on $\mcE_G$. We set $\mcE_i=\C[\mathfrak{t}]$, $ i\in I$.  
Let $w\in \W$
\[ w(\mcE_i)=\mcE_{iw^{-1}}, \quad \mcE_i=\C[\mathfrak{t}]\ni f \mapsto w\cdot f \in \C[\mathfrak{t}]=\mcE_{iw^{-1}} .\]
The isomorphism  $p:=\bigoplus_{i\in I} p_i$ defined by  
\[
\begin{aligned}
p_i\colon \C[\mathfrak{t}] &\to x_i^{-1} \C[\mathfrak{t}]\\
                                f&\mapsto x_i^{-1}(x_if)
\end{aligned} 
\]
gives the identification with the induced representation $\Ind_W^{\W}\C[\mathfrak{t}]$ which we described before. 

%
\end{rem}

\subsubsection{Calculations of some equivariant multiplicities}

In some situation one can actually say something on the images of algebraic cycle under the GKM-localization map, recall the 
\begin{satz} (multiplicity formular, \cite{Br}, section 3) 
Let $X$ equivariantly formal $T$-variety with a finite set of $T$-fixpoints $X^T$, by the localization theorem, 
\[
[X] =\sum_{x\in X^T} \La_x^X [\{x\}] \quad \in H_*^T(X)\otimes K
\]
where $\La_x^X\in K$. If $X$ is rationally smooth in $x$, then $\La_x^X\neq 0$ and $(\La_x^X)^{-1}= eu(X,x)\in H_{2n}^T(X)$, $\;$ 
$n=\dim_{\C}(X).$
\end{satz}

\begin{rem}
It holds for any $w\in \W$
\[
[\overline{Z^w}]= \sum_{i,j\in I} [\overline{Z_{i,j}^w}].
\]
Especially $1=[Z^e]=\sum_{i\in I} [Z_{i,i}^e]$ is the unit and $1_i=[Z_{i,i}^e]$ are idempotent elements, 
$1_i*1_j=0$ for $i\neq j$, $[Z_{i,j}]=1_i*[Z]*1_j$. 
In particular, for $s\in \S$ by lemma 11, we have 
\[
[\overline{Z^s}]= \sum_{i\in I\colon is=i} [\overline{Z_{i,i}^s}] + \sum_{i\in I\colon is \neq i} [{Z_{i,is}^s}].
\]

By the multiplicity formula we have 
\[ 
\begin{aligned}
{}[\overline{Z_{i,is}^s}] &= \begin{cases} 
\sum_{x\in W} \La_{xx_i, xx_is}^s\psi_{xx_i, xx_is}+ \La_{xx_i,xx_i}^s\psi_{xx_i,xx_i} &, \text{ if }i=is \\
\sum_{x\in W} \La_{xx_i,xx_is}^s \psi_{xx_i,xx_is}                                          &, \text{ if } is\neq i
\end{cases}\\
\text{with }& \La^s_{y,z}=(\eu(\overline{Z_{i,j}^s},\phi_{y,z}))^{-1}, \text{ for all }y,z\in \W\text{ as above} \\
[\overline{Z_{i,j}^w}] &= \begin{cases}
\sum_{x\in W} \La_{xx_i, xx_iw}^w\psi_{xx_i, xx_iw}+ \sum_{v < w}\La_{xx_i,xx_iv}^w\psi_{xx_i,xx_iv}  &,\text{ if } iw=j \\
0 &, \text{ if } iw\neq j
\end{cases}\\
\text{with }& \La^w_{xx_i,xx_iw} =(\eu(\overline{Z_{i,iw}^w}, \phi_{xx_i, xx_iw}))^{-1} \text{ for all }x\in W, 
\end{aligned}
\] 
\end{rem}

\subsubsection{Convolution on the fixed points}

The following key lemma on convolution products of $T$-fixed points
\begin{lemma}
For any $w, x,y\in \W$ it holds 
\[ \psi_{x,w} * \psi_w =\La_w \psi_x, \quad \psi_{x,w}*\psi_{w,y} = \La_{w} \psi_{x,y} \]
\end{lemma}

\paragraph{proof:}
We take $M_1=M_2=M_3=E$ and $Z_{1,2}:= \{\phi_{x,w} = ((0,xB),(0, wB)) \} \subset E\times E, Z_{2,3}:= \{ \phi_{w^\prime ,y} \} \subset E\times E $, then the 
set theoretic convolution gives 
\[ \{ \phi_{x,w} \} \circ \{ \phi_{w^\prime , y}\} = \begin{cases}\{\phi_{x,y}\}, &\text{ if } w=w^\prime \\
\emptyset, &\text{if }w\neq w^\prime
\end{cases}
\]
Similar, take $M_1=M_2=E, M_3=pt$, $Z_{12}:= \{\phi_{x,w} \}, Z_{23}={\phi_{w^\prime}}\times pt$, then 
\[\{ \phi_{x,w} \} \circ \{ \phi_{w}\}= \begin{cases} \{\phi_x\} &\text{if } w=w^\prime \\
\emptyset , &\text{ else} 
\end{cases}
\]
To see that we have to multiply with $\La_w$, we use the following proposition 
\begin{prop} (see \cite{CG}, Prop. 2.6.42, p.109)
Let $X_i\subset M, i=1,2$ be two closed (complex) submanifolds of a (complex) manifold with $X:=X_1\cap X_2$ is smooth and $T_xX_1\cap T_xX_2=T_xX$ for all $x\in X$. Then, we have 
\[ [X_1]\cap [X_2]= e (\mcT )\cdot [X] \] 
where $\mcT$ is the vector bundle $T_*M/(T_*X_1+T_*X_2)$ on $X$ and $e(\mcT )\in H^*(X)$ is the (non-equivariant) Euler class of this vector bundle, 
$\cap \colon H_*^{BM}(X_1)\times H_*^{BM}(X_2) \to H_*^{BM}(X)$ is the intersection pairing (cp. Appendix, or \cite{CG}, 2.6.15) and $\cdot$ on the right hand side stands for the $H^*(X)$-operation on the Borel-Moore homology (introduced in \cite{CG}, 2.6.40)     
\end{prop}  

Set $E_T:= E\times ^T ET, (\phi_x)_T:= \{\phi_x\}\times^T ET (\cong ET/T=BT)$.
We apply the proposition for $M=E_T^3$, $X_1:= (\phi_x)_T\times (\phi_w)_T \times E_T, X_2:= E_T \times (\phi_w)_T\times (\phi_y)_T$, $X_1\cap X_2\cong \{\phi_{x,y}\}_T (\cong BT)$, then 
$\mcT= (T_{\phi_w} E)\times^T ET$ and the (non-equivariant) Euler class is the top chern class of this bundle which is the $T$-equivariant top chern class 
of the constant bundle $T_{\phi_w} E$ on the point $\{\phi_{x,y}\}$. Since $T_{\phi_w} E =\bigoplus_{\la } \C_{\la}$ for one-dimensional $T$-representations 
$\C_{\la}$ with $t\cdot c:= \la (t) c,\;\; t\in T, c\in \C=\C_{\la}$. It holds 
\[
c_{top}^T(T_{\phi_w}E)= \prod_{\la} c_1^T(\C_{\la}) =\prod_{\la}\la = \La_w.
\]
Secondly, apply the proposition with $M= E_T^2\times (pt)_T, X_1=(\phi_x)_T\times (\phi_w)_T \times (pt)_T, X_2:=E_T\times (\phi_w)_T \times (pt)_T$, to see again 
$e(\mcT)=\La_w$.

\hfill $\Box$

Now we can give the missing proof of Corollary \ref{ConvOfZs}

\paragraph{proof of Corollary \ref{ConvOfZs}:}
By the lemma \ref{ConvolutionOfSubmodules} we know that there exists a $c\in \mcE_G$ such that $[\overline{Z^s_{i,j}}] * [\overline{Z^w_{j,k}}] =c\diamond [\overline{Z^{sw}_{i,k}}] $ 
in $ \mcZ_G^{\leq sw}/\mcZ_G^{< sw}$. We show that $c=1$. 
We pass with the forgetful map to $T$-equivariant Borel-Moore homology and tensor 
over $K={\rm Quot} (H_*^T(pt))$ and write $[\overline{Z^x_{s,t}}], x\in \W, s,t\in I$ for the image of the same named elements. 
Let $i,j,k\in I$ with $x_jwx_k^{-1} \in W$. 
\[
\begin{aligned}
{}[\overline{Z^s_{i,j}}] * [\overline{Z^w_{j,k}}] &= 
(\sum_{x\in W} \La_{xx_i, xx_is}^s\psi_{xx_i, xx_is}+ \La_{xx_i,xx_i}^s\psi_{xx_i,xx_i} )*\\   
& \;\; \; \; \; (\sum_{x\in W} \La_{xx_j, xx_jw}^w\psi_{xx_j, xx_jw}+ \sum_{v < w}\La_{xx_j,xx_jv}^w\psi_{xx_j,xx_jv} )\\
                                    &= \sum_{x\in W} \La_{xx_i, xx_is}^s \La_{xx_is,xx_isw}^w \La_{xx_is} \psi_{xx_i,xx_isw} + \underbrace{\cdots}_{\text{terms in }\mcZ_G^{< sw}}   
\end{aligned}
\]
Now, this has to be equal to $c\sum_{x\in W} \La_{xx_i,xx_isw}\psi_{xx_i,xx_isw}$ in $\mcZ_G^{\leq sw}/\mcZ_G^{< sw}$. 
Comparing coefficients at $x$ gives 
\[ 
\begin{aligned}
c&=
\frac{\eu (E_j,\phi_{xx_is}) \eu (\overline{Z_{i,k}^{sw}}, \phi_{xx_i, xx_isw})}{\eu (\overline{Z_{i,j}^s}, \phi_{xx_i,xx_is})\eu (\overline{Z_{j,k}^w}, \phi_{xx_is, xx_isw})}\\
&\\
&= \frac{\eu (\mathfrak{g}\cap {}^{xx_is}\mcU^- \oplus \mathfrak{g}\cap{}^{xx_i}\mcU^- \oplus\mathfrak{g}\cap{}^{xx_i}(\frac{\mcU^-}{\mcU^-\cap{}^{sw}\mcU^-}))}{\eu (\mathfrak{g}\cap{}^{xx_i}\mcU^- \oplus \mathfrak{g}\cap{}^{xx_i}(\frac{\mcU^-} {\mcU^-\cap{}^{s}\mcU^-}) \oplus\mathfrak{g}\cap{}^{xx_is}\mcU^- \oplus\mathfrak{g}\cap{}^{xx_i}(\frac{{}^{s}\mcU^-} {{}^{s}\mcU^-\cap{}^{sw}\mcU^-}))} \\
&\quad \cdot 
\prod_{l=1}^r \frac{\eu (V^{(l)} \cap {}^{xx_i}({}^s\mcU^{(l)}) \oplus V^{(l)} \cap {}^{xx_i}(\mcU^{(l)}\cap {}^{sw}\mcU^{(l)}))} 
{\eu (V^{(l)} \cap {}^{xx_i}(\mcU^{(l)}\cap{}^s\mcU^{(l)}) \oplus V^{(l)} \cap {}^{xx_i}({}^s\mcU^{(l)}\cap {}^{sw}\mcU^{(l)}))}  \\
&\\
&= \frac{ \eu ({}^x[\mathfrak{g}\cap{}^{x_i}(\frac{\mcU^-}{\mcU^-\cap{}^{sw}\mcU^-})])}{\eu ( {}^x[\mathfrak{g}\cap{}^{x_i}(\frac{\mcU^-} {\mcU^-\cap{}^{s}\mcU^-}) \oplus\mathfrak{g}\cap{}^{x_i}(\frac{{}^{s}\mcU^-} {{}^{s}\mcU^-\cap{}^{sw}\mcU^-})])}\\
&\cdot  \prod_{l=1}^r \frac{\eu ({}^x[V^{(l)} \cap {}^{x_i}(\frac{{}^s\mcU^{(l)}}{\mcU^{(l)} \cap {}^s\mcU^{(l)}}) \oplus V^{(l)} \cap {}^{x_i}(\frac{{}^{sw}\mcU^{(l)}}{{}^s\mcU^{(l)}\cap {}^{sw}\mcU^{(l)}})])}{\eu ({}^x[V^{(l)} \cap {}^{x_i}(\frac{{}^{sw}\mcU^{(l)}}{\mcU^{(l)}\cap{}^{sw}\mcU^{(l)}})] )} 
\end{aligned}
\] 
That for each $x$ and each $l\in \{1,\ldots , r\}$ the big two fraction in the product are equal to $1$ is a consequence of the following lemma.  
\hfill $\Box$

\begin{lemma} Let $T\subset\B\subset \G$ a maximal torus in a Borel subgroup in a reductive group (over$\C$), 
$F\subset Lie(\G)=\mcG$ a $\B$-subrepresentation. Let $(\W,\S)$ be the Weyl group for $(\G,T)$.   
Let $w\in \W, s\in \S$ such that $l(sw)=l(w)+1$, then it holds for any $x\in \W$  
\[ {}^x(\frac{{}^sF}{F\cap {}^sF} \oplus {}^s( \frac{{}^wF}{F\cap {}^wF} )) \cong {}^x(\frac{{}^{sw}F}{F\cap {}^{sw}F}). \]
In particular, this holds also for $F= \mathfrak{u}^{-}$. 
\end{lemma}

\paragraph{proof:}
Let $\Phi_F:= \{ \al \in \Hom(\mathfrak{t}, \C) \mid F_{\al}\neq 0\} \subset \underline{\Phi}$, $\Phi^+(y):= \underline{\Phi}^+\cap y(\underline{\Phi}^-)$, $\Phi_F^+(y):=\Phi_F\cap  \Phi^+(y), y\in \W$ where $\underline{\Phi},\underline{\Phi}^+,\underline{\Phi}^-$ are the set of roots (of $T$ on $\mcG$), positive roots, negative roots respectively.\\ 
The assumption $l(sw)=l(w)+ 1$ implies $\Phi^+_F(sw)= s\Phi^+_F(w) \sqcup \Phi^+_F(s)$ and 
for $\Phi_F^{-}(y):=-\Phi_F^+(y)$,  $\Phi_F(y):= \Phi_F^+(y)\cup \Phi^-_F(y)= \Phi_F\setminus (\Phi_F \cap y\Phi_F)$ it holds 
$ \Phi_F(sw)= s\Phi_F(w) \sqcup \Phi_F(s)$ and 
for any $x\in W$ it holds $ x\Phi_F(sw)= x(s\Phi_F(w) \sqcup \Phi_F(s))$. Now, the weights of ${}^x(\frac{{}^{sw}F}{F\cap {}^{sw}F})$ are $x\Phi_F(sw)$, 
the weights of ${}^x(\frac{{}^sF}{F\cap {}^sF} \oplus {}^s( \frac{{}^wF}{F\cap {}^wF} ))$ are $x(s\Phi_F(w) \sqcup \Phi_F(s)) $.
\hfill $\Box$

\subsection{Generators for $\mcZ_G$}

Let $J=\emptyset $. Recall, we denote the right $\W$-operation on $I= W \setminus \W$ by $(i,w)\mapsto iw$, $i\in I, w\in \W$.\\
For $i\in I$ we set $\mcE_i:=  H_G^*(E_i)=\C[\mathfrak{t}]=\C[x_i(1),\ldots , x_i(m)]$, we write 
\[w(\al_s)=w(\al_s(x_{iw^{-1}}(1), \ldots , x_{iw^{-1}}(m)))\in \mcE_{iw^{-1}}\] 
for the element corresponding to the root $w(\al_s), s\in \S, w\in \W$ without mentioning that it depends on $i\in I$. \\
We define a collection of elements in $\mcZ_G$
\[
\begin{aligned}
1_i &:= [ Z_{i,i}^e ] \\
z_{i}(t) &:= x_i(t)  \in \mcZ_G^{\leq e} (\subset \mcZ_G )\\
\si_{i}(s) &:=  [\overline{Z^s_{i,j}}] \in \mcZ_G^{\leq s}, \text{ where }is=j  
\end{aligned}
\]
where we use that $\mcE_i \subset \mcZ_G^{\leq e} \subset \mcZ_G$ and the degree of $x_i(t)$ is $2$ in $H_{[*]}^G(Z)$, 
see Lemm 6 and the definition of the grading (just before theorem 2.1) . It is also easy to 
see that $1_i\in H_{[0]}^G(Z)$ because $\deg 1_i = 2e_i - 2\dim Z_{i,i}^e = 0$. Furthermore, the degree of $\si_i(s)$ is  
\[ e_{is}+ e_i - 2\dim Z_{i, is}^s= \begin{cases} 2 \deg q_i(s) -2 , & \text{ if } is=i \\
                                                  2\deg q_i(s) , &\text { if } is \neq i. \end{cases}.\]
Recall $\mcZ_G\inj \End (\mcE_G )=\End (\bigoplus_{i\in I} \mcE_i)$ from \cite{VV2}, remark after Prop.3.1, p.12. Let us denote by $\widetilde{1_i},\widetilde{z_{i}(t)},\widetilde{\sigma_{i} (s)}$ be the images of $1_i,z_{i}(t),\sigma_{i} (s)$.\\

\begin{prop} 
Let $k\in I$, $f\in \mcE_k $, $\al_s \in \underline{\Phi}^+$ be the positive root such that $s(\al_s)=-\al_s$. 
It holds 
\[
\begin{aligned}
\widetilde{1_i} (f) &:= 1_i*f = \begin{cases} f,  &\text{ if }i=k, \\
                                0, &\text{ else. }
\end{cases} \\
\widetilde{z_{i}(t)}(f)&:=z_{i}(t)*f = \begin{cases} x_i(t)f ,  &\text{ if }i=k, \\
                       0, & \text{ else. }
\end{cases} \\
\widetilde{\sigma_{i} (s)} (f) &:= \begin{cases} q_i(s)\frac{s(f)-f}{\al_s} , \quad  &\text{ if }i=is=k , \\
                                                    q_{i}(s) s(f) &\text{ if }i\neq is =k, \\
                       0, & \text{ else. }
\end{cases}     \\
\text{ for }& \mcU=\Lie (\U)^{\oplus r} \text{ this looks like}\\
\widetilde{\sigma_{i} (s)} (f) &:= \begin{cases} \al_s^{h_i(s)}\frac{s(f)-f}{\al_s} , \quad  &\text{ if }i=is=k, \\
                                                    \al_s^{h_{i}(s)} s(f) &\text{ if }i\neq is =k, \\
                       0, & \text{ else. }
\end{cases}                     
\end{aligned}
\]
We write $\delta_s:=\frac{s-1}{\al_s}$, it is the BGG-operator from \cite{De}, i.e. 
for $is=i, f\in \mcE_i$,  \\ $\si_i(s) (f) = q_i(s) \delta_s(f)$. 
\end{prop}

\paragraph{proof:} 
Consider the following two maps 
\[
\begin{aligned}
\Theta \colon \mcE_G\to \mcE_T \to \mcE_T \otimes K &\to \bigoplus_{w\in \W} K\psi_w \\
\mcE_k\ni \;\; f  & \mapsto \sum_{w\in Wx_k} w(f) \La_w^{-1} \psi_w 
\end{aligned}
\]
\[
\begin{aligned}
C\colon &\bigoplus_{w\in \W} K\psi_w \to \bigoplus_{w\in \W} K\psi_w\\
               \psi_w &\mapsto [\overline{Z^s_{i,is}}]*\psi_w =
               \begin{cases} 
               (\sum_{x\in W}\La_{xx_i,xx_i}^s \psi_{xx_i,xx_i} + \La_{xx_i,x{}^{x_i}sx_i}^s\psi_{xx_i,x{}^{x_i}sx_i}) *\psi_w   & \\ 
                      \quad \quad \quad = \La_{w,w}^s \La_w \psi_w + \La_{ws,w}^s \La_w\psi_{ws},  & \text{ if }w\in Wx_i, \; i=is\\
               (\sum_{x\in W}\La_{xx_i,x{}^{x_i}sx_i}^s\psi_{xx_i,x{}^{x_i}sx_i}) *\psi_w   & \\ 
                      \quad \quad \quad = \La_{ws,w}^s \La_w\psi_{ws},  &\text{ if } w\in Wx_is,\; i\neq is \\
                   0 , &\text{ if } w\notin Wx_is  
               \end{cases}                   
\end{aligned}
\]
To calculate $[\overline{Z^s_{i,is}}]*f, f\in \mcE_k$ it is enough to calculate $[\overline{Z^s_{i,is}}]*\Theta(f)= C(\Theta(f))$ because $\Theta $ is an injective algebra homomorphism. 
\[
C\Theta (f)= \begin{cases} \delta_{is,k}\sum_{w\in Wx_i} [w(f) \La_{w,w}^s +  w(sf)\La_{w,ws}^s ]\psi_w, & \text{ if }i=is \\
                      \delta_{is,k} \sum_{w\in Wx_i} [w(sf)\La_{w,ws}^s ]\psi_w, & \text{ if }i\neq is 
       \end{cases}                
\] 
Now, recall,
\begin{itemize}
\item[(1)]  If $i=is=k$
\[
\begin{aligned}
C\Theta(f) &= \sum_{w\in Wx_i} w[ q_i(s)\frac{s(f)-f}{\al_s} ] \La_w^{-1} \psi_w \\
      &= \Theta(q_i(s)\frac{s(f)-f}{\al_s} ) 
\end{aligned}
\]
Once we identify $\mcE_k=\C[\mathfrak{t}], \; k\in I$, we see that $\si_i(s)\colon \mcE_G\to \mcE_G$ is the zero map on the $k$-th summand, $k\neq i$ and on the $i$-th summand 
\[ 
\begin{aligned}
\C[\mathfrak{t}] &\to \C[\mathfrak{t}] \\
f &\mapsto q_i(s) \frac{s(f)-f}{\al_s} 
\end{aligned}
\]  
\item[(2)] If $i\neq is=k$,
\[
\begin{aligned}
C\Theta (f)&= \sum_{w\in Wx_i} [w(sf)\La_{w,ws}^s ]\psi_w \\
     & = \Theta(q_i(s) s(f))
\end{aligned}
\]
Once we identify 
$\mcE_k=\C[\mathfrak{t}]$, we see that $\si_i(s)\colon \mcE_G\to \mcE_G$ is the zero map on the $k$-th summand, $k\neq is$ and on the $is$-th summand 
it is the map 
\[
\begin{aligned}
\C[\mathfrak{t}] &\to \C[\mathfrak{t}]\\
f&\mapsto q_i(s) s(f)
\end{aligned}
\]

\end{itemize}
\hfill $\Box $

\begin{lemma}
The algebra $\mcZ_G$ is generated as $\La_G$-algebra by the elements 
\[1_i, i\in I, \quad z_i(t), 1\leq t\leq rk(T), i\in I,\quad   \sigma_i(s), s\in \S, i\in I .\] 
\end{lemma}

\paragraph{proof:}
It follows from the cellular fibration property that $\mcZ_G$ is generated by $1_i, i\in I, \quad z_i(t), 1\leq t\leq rk(T), i\in I, [\overline{Z_{i,j}^w}], w\in \W$. By corollary \ref{ConvOfZs} it follows that one can restrict to the case $w\in \S$, more precisely as free $H_G^*(E)$-module it can be generated by 
\[ \si (w):= \si(s_1)*\cdots \si(s_t), w\in \W, w=s_1\cdots s_t \text{ reduced expression }, \si(s):= \sum_{i\in I }\si_i(s),\]
and this basis has a unitriangular base cange to the basis given by the $[\overline{Z^w}]$. 

\hfill $\Box$

Furthermore, we consider \[\Phi\colon \bigoplus_{i\in I} \C[x_i(1),\ldots x_i(n)]\cong \bigoplus_{i\in I} \C[z_i(1),\ldots z_i(n)], \; x_i(t)\mapsto z_i(t)\] 
as the left $\W$-module $\Ind_W^{\W} \C[\mathfrak{t}]$, we fix the polynomials 
\[ c_{i}(s,t):= \Phi(\si_i(s)(x_i(t))) \;\in \bigoplus_{i\in I} \C[z_i(1),\ldots z_i(n)], \quad i\in I, \; 1\leq t\leq n,\; s\in \S. \]

\begin{prop} \label{relations}
Under the following assumption for the data $(\G, \B , \mcU= (\Lie (\U))^{\oplus r}, H,V), J=\emptyset\colon $
Let $\S\subset \W=Weyl(\G,T)$ be the simple reflections, we assume for any $s,t\in \S$ 
\begin{itemize}
\item[(B2)] If the root system spanned by $\al_s, \al_t$ is of type $B_2$ (i.e. $stst=tsts$ is the minimal relation), then for every $i\in I$ such that 
$is=i=it$ it holds $h_i(s), h_i(t)\in \{0,1,2\}$. 
\item[(G2)] If the root system spanned by $\al_s, \al_t$ is of type $G_2$ (i.e. $ststst=tststs$ is the minimal relation), then for every $i\in I$ such that 
$is=i=it$ it holds $h_i(s)=0=h_i(t)$. 
\end{itemize}  
Then the generalized quiver Hecke algebra for $(\G, \B , \mcU= (\Lie (\U))^{\oplus r}, H,V)$ is the $\C$-algebra with generators 
 \[1_i, i\in I, \quad  z_i(t), 1\leq t\leq n=rk(T), i\in I, \quad  \si_i(s), s\in \S, i\in I\] 
 and relations 
\begin{itemize}

\item[(1)] (\emph{orthogonal idempotents})
\[
\begin{aligned}
1_i 1_j  &=\delta_{i,j} 1_i, \\
1_i z_i(t) 1_i  &= z_i(t), \\
1_i \si_i(s) 1_{is} &= \si_i(s)
\end{aligned}
\]
\item[(2)] (\emph{polynomial subalgebras}) 
\[z_i(t) z_i(t^\prime)= z_i(t^\prime) z_i(t)\]

\item[(3)] (\emph{ relation implied by $s^2=1$}) 
\[\si_i(s) \si_{is}(s) = \begin{cases}  0 &, \text{ if }is=i,\;  h_i(s) \text{ is even } \\
                                    -2\al_s^{h_i(s)-1} \si_i(s) & , \text{ if }is=i,\;  h_i(s) \text{ is odd } \\
                                    (-1)^{h_{is}(s)} \al_s^{h_i(s)+h_{is}(s)} &, \text{ if } is\neq i
\end{cases}
\]

\item[(4)] (\emph{straightening rule})\\
\[ 
\si_i(s) z_i(t) - s(z_i(t))\si_i(s) = \begin{cases} c_{i}(s,t), &\text{, if } is=i \\
0 &\text{, if } is\neq i. \end{cases}
\]

\item[(5)] (\emph{braid relations})\\
Let $s,t\in \S, st=ts$, then 
\[  \si_i(s)\si_{is}(t) =  \si_i (t) \si_{it}(s)
\]
 Let $s, t\in \S$ not commuting such that $x:=sts \cdots= tst\cdots $ minimally, $i\in I$. 
There exists explicit polynomials $(Q_w)_{w< x}$ in $\al_s, \al_t\in \C[\mathfrak{t}]$ such that 
\[ 
\si_i (sts \cdots ) -\si_i (tst\cdots ) = \sum_{w<x} Q_w \si_i (w) 
\]
(observe that for $w<x$ there exists just one reduced expression).  

\end{itemize}

\end{prop}

\paragraph{proof:} For the convenience of the reader who wants to check the relations for the generators of $\mcZ_G$, 
we include the detailed calculations. (1), (2) are clear. Let always $f\in \C[\mathfrak{t}]\cong \mcE_{is}$.  We will 
use as shortage $\delta_s (f):= \frac{s(f)-f}{\al_s}$ and use that these satisfy the usual relations of BGG-operators (cp. \cite{De}).  
\begin{itemize}

\item[(3)]
If $is=i$, then
\[ \begin{aligned} \si_i (s)\si_i(s)(f) & = \al_s^{h_i(s)} \delta_s (\al_s^{h_i(s)} \delta_s(f)) \\ &= \al_s^{h_i(s)} \delta_s (\al_s^{h_i(s)}) \delta_s(f) =  [(-1)^{h_i(s)}-1]\al_s^{h_i(s)-1} \si_i(s)(f).
\end{aligned}
\]

If $is\neq i$, then 
\[ \si_i(s)(f)\si_{is}(s) = \al_s^{h_i(s)} s(\al_s^{h_{is}(s)}) s(s(f)) =  (-1)^{h_{is}(s)} \al_s^{h_i(s)+h_{is}(s)} f. \] 

\item[(4)]
(straightening rule)

The case $is\neq i$ is clear by definition. Let $is=i$, then the relation follows directly from the product rule for BGG-operators, which 
states $\delta_s( xf) ) = \delta_s(x) f+s(x) \delta_s(f),\;\; x, f\in \C[\mathfrak{t}].$

\item[(5)] (braid relations) \\
$s,t\in \S, st=ts$, $f\in \C[\mathfrak{t}]$, to prove  
\[  \si_i(s)\si_{is}(t) (f) =  \si_i (t) \si_{it}(s)(f) 
\]
we have to consider the following four cases. We use the following: $t(\al_s)=\al_s, s(\al_t)=\al_t, h_i(s)=h_{it}(s), h_{i}(t)=h_{is}(t), 
\delta_s(\al_t^{h_i(t)})=0=\delta_t(\al_s^{h_i(s)})$.  
\begin{itemize}
\item[1.] $is=i, it=i$,  use $\delta_s \delta_t =\delta_t \delta_s$
\[
\begin{aligned}
 \si_i(t)\si_i(s)(f)=\al_t^{h_{i}(t)} \delta_t(\al_s^{h_i(s)} \delta_s(f)) &= \al_s^{h_i(s)}\al_t^{h_i(t)} \delta_s\delta_t (f) \\
 &= \al_s^{h_i(s)}\al_t^{h_i(t)} \delta_t\delta_s (f) \\
 &= \al_s^{h_{i}(s)} \delta_s(\al_t^{h_i(t)} \delta_t(f)) =\si_i(s)\si_i(t) (f) 
 \end{aligned} .
 \]
\item[2.] $is=i, it\neq i$, use $\delta_s t =t\delta_s$\\
\[
\begin{aligned}
\si_i(t) \si_{it}(s) (f) = \al_t^{h_i(t)} t(\al_s^{h_{it}(s)}\delta_s(f) ) &=  \al_t^{h_i(t)} \al_s^{h_i(s)} t\delta_s(f) \\
  &=  \al_t^{h_{is}(t)} \al_s^{h_{i}(s)} \delta_s(t(f))  \\
  &= \al_s^{h_{i}(s)}\delta_s(\al_t^{h_{is}(t)} t(f)) =\si_i(s) \si_i (t) (f) 
\end{aligned}.
\] 
\item[3.] $is\neq i, it=i$, follows by symmetry from the last case.
\item[4.] $is\neq i, it\neq i$. \\
\[
\begin{aligned}
\si_i(t)\si_{it}(s)(f) =
\al_t^{h_i(t)} t(\al_s^{h_{it}(s)} s(f)) = \al_s^{h_i(s)} s(\al_t^{h_{is}(t)} t(f))\\
= \si_i(s) \si_{is}(t)(f). 
\end{aligned}
\]
\end{itemize}

\noindent
Let $st\neq ts$. There are three different possibilties, either 
\begin{itemize}
\item[(A)] $sts=tst$  \hfill (type $A_2$)
\item[(B)] $stst=tsts$  \hfill (type $B_2$) 
\item[(C)] $ststst=tststs$ \hfill (type $G_2$)
\end{itemize}
We write $Stab_i:=\{ w\in \langle s,t\rangle \mid iw=i\}$. For each case we go through the subgrouplattice to calculate explicitly the polynomials $Q_w$. 

\paragraph{(A) $sts=tst\colon $} $\langle s,t\rangle \cong S_3$, $s(\al_t)=t(\al_s)=\al_s+\al_t$. 
We have five (up to symmetry between $s$ and $t$) subgroups to consider. Always, it holds 
\[
h_{is}(t)=h_{it}(s), h_{ist}(s)=h_i(t), h_{its}(t)=h_i(s)
\]
which implies an equality which we use in all five cases 
\[ 
\begin{aligned}
 \al_s^{h_i(s)}s(\al_t^{h_{is}(t)})st(\al_s^{h_{ist}(s)}) &= \al_s^{h_i(s)}(\al_s+\al_t)^{h_{it}(s)}\al_t^{h_{ist}(s)} \\
 &= \al_t^{h_i(t)} t(\al_s^{h_{it}(s)}) ts(\al_t^{h_{its}(t)}) 
\end{aligned}
\] 
\begin{itemize}
\item[A1.] $Stab_i= \langle s,t\rangle $, this implies $h_i(s)=h_i(t)=:h$ by definition ($x_i(\al_s)\in \Phi_{V^{(k)}}$ if and 
only if $ x_{it}(\al_s)=x_i(\al_s+\al_t)=x_{is}(\al_t) \in \Phi_{V^{(k)}} \Leftrightarrow x_i(\al_t)\in \Phi_{V^{(k)}}$) and as a consequence we get \\
$\al^h\delta_s(\al_t^h t(\al_s^h))=0$. This simplifies the equation to 
\[ 
\si_i(s)\si_i(t)\si_i(s)-\si_i(t)\si_i(s)\si_i(t) = \delta_s(\al_t^h\delta_t(\al_s^h)) \si_i(s) - \delta_t(\al_s^h\delta_s(\al_t^h)) \si_i(t)
\]
note that $Q_s:= \delta_s(\al_t^h\delta_t(\al_s^h)) , Q_t:=-\delta_t(\al_s^h\delta_s(\al_t^h))$ are polynomials in $\al_s, \al_t$.
\item[A2.] $Stab_i=\langle s \rangle $ (analogue $Stab_i=\langle t\rangle $). It holds $itst=its$. We use in this case \\
 $h_{is}(t)=h_i(t), h_{ist}(s)=h_{it}(s)= h_i(t), h_{its}(t) =h_i(s)$. 
\[
\begin{aligned}
\si_{i}(s) \si_i(t)\si_{it}(s)(f) &- \si_{i}(t)\si_{it}(s) \si_{its}(t)(f) \\
&=  \al_s^{h_{i}(s)}\delta_s(\al_t^{h_{is}(t)}t(\al_s^{h_{ist}(s)}) ts (f)) - \al_t^{h_i(t)} t(\al_s^{h_{it}(s)}) ts(\al_t^{h_{its}(t)}) ts \delta_t(f)\\ 
&= \al_s^{h_i(s)} \delta_s(\al_t^{h_{is}(t)}t(\al_s^{h_{ist}(s)})) ts(f) + \al_s^{h_i(s)} s(\al_t^{h_{is}(t)}) st(\al_s^{h_{ist}(s)}) \delta_s (ts (f)) \\
&=  \quad -   \al_t^{h_i(t)} t(\al_s^{h_{it}(s)}) ts(\al_t^{h_{its}(t)}) ts \delta_t(f) =0
\end{aligned}
\]
Since $st\delta_s=\delta_t st$ and $ \delta_s (\al_t^{h} t(\al_s)^{h})=0$.

\item[A3.] $Stab_i=\langle sts\rangle $, then $ist=is, its=it$. 
\[
\begin{aligned}
\si_{i}(s) \si_{is}(t)\si_{is}(s)(f) &- \si_{i}(t)\si_{it}(s) \si_{it}(t)(f) \\
&=  \al_s^{h_{i}(s)} s(\al_t^{h_{is}(t)} \delta_t(\al_s^{h_{is}(s)}s(f))) - \al_t^{h_{is}(s)} t(\al_s^{h_{is}(t)} \delta_s(\al_t^{h_{i}(s)}t(f)))  \\
&= [\al_s^{h_{i}(s)} s(\al_t^{h_{is}(t)}) s(\delta_t(\al_s^{h_{is}(s)})) - \al_t^{h_{is}(s)} t(\al_s^{h_{is}(t)}) t(\delta_s(\al_t^{h_{i}(s)}))]\cdot f
\end{aligned}
\]
using $t\delta_s t= s\delta_t s$.
\item[A4.] $Stab_i=\{1\}$ (and the same for $Stab_i=\langle st\rangle $)
\[
\begin{aligned}
\si_{i}(s) \si_{is}(t)\si_{ist}(s) &- \si_{i}(t)\si_{it}(s) \si_{its}(t) \\
&=  \al_s^{h_{i}(s)}s(\al_t^{h_{is}(t)})st(\al_s^{h_{ist}(s)}) sts - \al_t^{h_{i}(t)} t(\al_s^{h_{it}(s)}) ts(\al_t^{h_{its}(t)}) tst\\
&= 0 
\end{aligned}
\]
\end{itemize} 
\paragraph{(B) $stst=tsts\colon $ } $<s,t>\cong D_4 $(order is $8$),  
\[
\begin{aligned}
t(\al_s)=\al_s+\al_t, &\quad st(\al_s)=\al_s+\al_t, \quad tst(\al_s)=\al_s   \\
s(\al_t)= 2\al_s +\al_t, &\quad  ts(\al_t) =2\al_s+\al_t,\quad  sts(\al_t) =\al_t.
\end{aligned}
\] 
Here we have to consider ten different cases because $D_4$ has ten subgroups. It always holds the following 
\[ h_{itst}(s)=h_i(s), h_{its}(t)=h_{is}(t), h_{it}(s)=h_{ist}(s), h_{ists}(t)=h_i(t) \]
which implies 
\[ 
\al_s^{h_{i}(s)}s(\al_t)^{h_{is}(t)} st (\al_s)^{h_{ist}(s)} sts(\al_t)^{h_{ists}(t)} = \al_t^{h_{i}(t)}t(\al_s)^{h_{it}(s)} ts (\al_t)^{h_{its}(t)} tst(\al_s)^{h_{itst}(s)}
\]
This will be used in all cases, it is particular easy to see that for 
\[Stab_i=\{1\},\quad  Stab_i=\{1,ts,st,stst\}, \quad Stab_i=\{ 1, stst \}\] 
we obtain that the difference is zero from the above equality. 
Let us investigate the other cases. Furthermore, the following is useful to notice
\[
\delta_s(t(\al_s)^h)=0, \quad \delta_t(s(\al_t)^h)=0 
\]
\begin{itemize}
\item[B1.] $Stab_i=\langle s,t\rangle $. We prove the following 
\[ 
\begin{aligned}
\si_i(s)\si_i(t)\si_i(s)\si_i(t) (f) &= Q_{st}\sigma_i(s)\sigma_i(t) + \al_s^{h_{i}(s)}s(\al_t)^{h_{i}(t)} st (\al_s)^{h_{i}(s)} sts(\al_t)^{h_i(t)} \delta_{stst}(f)\\
& \; \; +  \al_s^{h_i(s)}s(\al_t^{h_i(t)}) st(\al_s^{h_i(s)}) \delta_{sts}(\al_t^{h_i(t)}) \delta_t
\end{aligned}
\]
with $Q_{st}= \delta_s(\al_t^{h_i(t)}) \delta_t(\al_s^{h_i(s)}) + s(\al_t^{h_i(t)}) \delta_{st}(\al_s^{h_i(s)}) + t(\al_s^{h_i(s)}) \delta_{ts}(\al_t^{h_i(t)}) =Q_{ts}$ is a polynomial in $\al_s,\al_t$.
By a long direct calculation (applying the product rule for the $\delta_s$) several times 
\[ 
\begin{aligned}
\si_i(s)\si_i(t)&\si_i(s)\si_i(t) (f) = \al_s^{h_i(s)} \delta_s(\al_t^{h_i(t)} \delta_t (\al_s^{h_i(s)} \delta_s(\al_t^{h_i(t)}))) \delta_t(f)  \\
&+[\al_s^{h_i(s)} s(\al_t^{h_i(t)}) s\delta_t(\al_s^{h_i(s)}\delta_s(\al_t^{h_i(t)})) + \al_s^{h_i(s)}\delta_s(\al_t^{h_i(t)}\delta_t(\al_s^{h_i(s)}s(\al_t^{h_i(t)}))) ] \delta_{st}(f) \\
&+ \al_s^{h_i(s)}\delta_s(\al_t^{h_i(t)} t(\al_s^{h_i(s)})ts(\al_t^{h_i(t)})) \delta_{tst}(f) \\ 
&+ \al_s^{h_{i}(s)}s(\al_t)^{h_{i}(t)} st (\al_s)^{h_{i}(s)} sts(\al_t)^{h_i(t)} \delta_{stst}(f)
\end{aligned}
\]
We have a look at the polynomials occurring in front of the $\delta_w$:
\begin{itemize}
\item[$w=t\colon $] by the product rule 
\[ 
\begin{aligned}
\al_s^{h_i(s)}& \delta_s(\al_t^{h_i(t)} \delta_t (\al_s^{h_i(s)} \delta_s(\al_t^{h_i(t)}))) = \al_s^{h_i(s)} \delta_s(\al_t^{h_i(t)})^2\delta_t(\al_s^{h_i(s)})\\
&+\al_s^{h_i(s)}s(\al_t^{h_i(t)}) \delta_s(\al_t^{h_i(t)}) \delta_{st}(\al_s^{h_i(s)}) + \al_s^{h_i(s)} t(\al_s^{h_i(s)}) \delta_s(\al_t^{h_i(t)}) \delta_{ts}(\al_t^{h_i(t)})\\
&+ \al_s^{h_i(s)}s(\al_t^{h_i(t)}) st(\al_s^{h_i(s)}) \delta_{sts}(\al_t^{h_i(t)}) 
\end{aligned}
\]
\item[$w=st\colon $] 
\[ 
\begin{aligned}
\al_s^{h_i(s)}& s(\al_t^{h_i(t)}) s\delta_t(\al_s^{h_i(s)}\delta_s(\al_t^{h_i(t)})) + \al_s^{h_i(s)}\delta_s(\al_t^{h_i(t)}\delta_t(\al_s^{h_i(s)}s(\al_t^{h_i(t)}))) = \\
\al_s^{h_i(s)}& s(\al_t^{h_i(t)})s\delta_s(\al_t^{h_i(t)})s\delta_t(\al_s^{h_i(s)}) + \al_s^{h_i(s)} s(\al_t^{h_i(t)})st(\al_s^{h_i(s)}) \delta_{ts}(\al_t^{h_i(t)}) \\
+ \al_s^{h_i(s)}& s(\al_t^{h_i(t)})\delta_s(\al_t^{h_i(t)})\delta_t(\al_s^{h_i(s)}) + \al_s^{h_i(s)}s(\al_t^{h_i(t)}) \delta_s(s(\al_t^{h_i(t)})) s\delta_t(\al_s^{h_i(s)}) \\
+ \al_s^{h_i(s)}& s(\al_t^{h_i(t)})^2\delta_{st}(\al_s^{h_i(s)}) \\
= \al_s^{h_i(s)}& s(\al_t^{h_i(t)}) [ \delta_s(\al_t^{h_i(t)}) \delta_t(\al_s^{h_i(s)}) + s(\al_t^{h_i(t)}) \delta_{st}(\al_s^{h_i(s)}) + t(\al_s^{h_i(s)}) \delta_{ts}(\al_t^{h_i(t)}) ] \\
=\al_s^{h_i(s)}& s(\al_t^{h_i(t)}) Q_{st}
\end{aligned}
\]
using $s(\delta_s(\al_t^h))=\delta_s(\al_t^h)$ and $\delta_s(s(\al_t)^h)= -\delta_s(\al_t^h)$. 
\item[$w=tst\colon $]
\[
\al_s^{h_i(s)}\delta_s(\al_t^{h_i(t)} t(\al_s^{h_i(s)})ts(\al_t^{h_i(t)}))=0 
\]
\end{itemize}
Now, look at $\si_i(s)\si_i(t) (f) = \al_s^{h_i(s)}\delta_s(\al_t^{h_i(t)})\delta_t(f) +\al_s^{h_i(s)} s(\al_t^{h_i(t)}) \delta_{st}(f)$,
which implies 
\[\al_s^{h_i(s)} s(\al_t^{h_i(t)}) Q_{st}\delta_{st}(f) = Q_{st}\si_i(s)\si_i(t)(f) - \al_s^{h_i(s)} \delta_s(\al_t^{h_i(t)}) Q_{st} \delta_t(f)\]
replace the previous expression and compare coefficients in front of $\delta_t(f)$ again gives the polynomial  
\[ 
\begin{aligned}
\al_s^{h_i(s)}& \delta_s(\al_t^{h_i(t)} \delta_t (\al_s^{h_i(s)} \delta_s(\al_t^{h_i(t)}))) - \al_s^{h_i(s)} \delta_s(\al_t^{h_i(t)}) Q_{st} = \\
\al_s^{h_i(s)}& s(\al_t^{h_i(t)}) st( \al_s^{h_i(s)}) \delta_{sts}(\al_t^{h_i(t)})
\end{aligned}
\]
We conclude 
\[
\begin{aligned}
\si_i(s)&\si_i(t)\si_i(s)\si_i(t) -\si_i(t)\si_i(s)\si_i(t)\si_i(s) = Q_{st}\si_i(s)\si_i(t) - Q_{st} \si_i(t)\si_i(s) \\
&+ \al_s^{h_i(s)} s(\al_t^{h_i(t)}) st( \al_s^{h_i(s)}) \delta_{sts}(\al_t^{h_i(t)})\delta_t - \al_t^{h_i(t)} t(\al_s^{h_i(s)}) ts( \al_t^{h_i(t)}) \delta_{tst}(\al_s^{h_i(s)})\delta_s 
\end{aligned}
\]
Since $\delta_{sts}(\al_t^h) =0= \delta_{tst}(\al_s ^k)$ for $h,k\in \{0,1,2\}$ since the maps $\delta_{sts},\delta_{tst}$ map polynomials of degree $d$ to polynomials of degree $d-3$ or to zero, the claim follows. In general, if we localize to $\C[\mathfrak{t}][\al_t^{-1}, \al_s^{-1}]$ we could still have the analogue statement.   

\item[B2.] $Stab_i=\langle s\rangle $ (analogue $Stab_i=\langle t\rangle $) and use $\delta_s(\al_t^{h_{is}(t)} t(\al_s^{h_{ist}(s)}) ts(\al_t^{h_{ists}(t)})) =0$ to see 
\[
\begin{aligned}
&\si_{i}(s)\si_{is}(t)\si_{ist}(s)\si_{ists}(t)- \si_{i}(t)\si_{it}(s)\si_{its} (t)\si_{itst}(s) \\
&= \al_s^{h_{i}(s)} \delta_s(\al_t^{h_{is}(t)} t(\al_s^{h_{ist}(s)}) ts(\al_t^{h_{ists}(t)}) tst (f)) \\
&\quad  - \al_t^{h_{i}(t)}t(\al_s^{h_{it}(s)}) ts(\al_t^{h_{its}(t)})tst(\al_s^{h_{itst}(s)}) tst\delta_s(f)\\
&=\al_s^{h_{i}(s)} s(\al_t^{h_{is}(t)}) st(\al_s^{h_{ist}(s)}) sts(\al_t^{h_{ists}(t)}) \delta_s tst (f)) \\
&\quad - \al_t^{h_{i}(t)}t(\al_s^{h_{it}(s)}) ts(\al_t^{h_{its}(t)})tst(\al_s^{h_{itst}(s)}) tst\delta_s(f) \\
&=0
\end{aligned}
\]

because $tst\delta_s=\delta_s tst$. 

\item[B3.] $Stab_i=\{1, sts\}$ (analogue $Stab_i=\{1, tst\}$). It holds $its=itst, is=ist$. We have 
\[
\begin{aligned}
{}[\si_{i}(s)&\si_{is}(t) \si_{ist}(s) \si_{ists}(t) - \si_i(t)\si_{it}(s) \si_{its}(t)\si_{itst}(s)](f)  \\
=&\al_s^{h_{i}(s)} s(\al_t^{h_{is}(t)}) s\delta_t(\al_s^{h_{ist}(s)} s(\al_t^{h_{ists}(t)}) st(f)) \\
&\quad - \al_t^{h_i(t)}t(\al_s^{h_{it}(s)}) ts(\al_t^{h_{its}(t)}) ts\delta_t(\al_s^{h_{itst}(s)}s(f))  \\
=&\al_s^{h_{i}(s)} s(\al_t^{h_{is}(t)})[ s\delta_t(\al_s^{h_{ist}(s)} s(\al_t^{h_{ists}(t)})) t(f) + st(\al_s^{h_{ist}(s)}) sts (\al_t^{h_{ists}(t)}) s\delta_t (st(f))] \\
&\quad -  \al_t^{h_i(t)} t(\al_s^{h_{it}(s)})ts(\al_t^{h_{its}(t)}) [ts\delta_t(\al_s^{h_{itst}(s)}) t(f) +   tst(\al_s^{h_{itst}(s)}) ts\delta_t(s(f))] \\
=&[\al_s^{h_{i}(s)} s(\al_t^{h_{is}(t)}) s\delta_t(\al_s^{h_{ist}(s)}) -t(\al_s^{h_{it}(s)}) ts(\al_t^{h_{its}(t)}) ts \delta_t(\al_s^{h_{itst}(s)})] \si_i(t) (f)
\end{aligned}
\]
using $s\delta_t st = ts\delta_ts$ and $s\delta_t(\al_s^{h_{ist}(s)}s(\al_t^{h_{ists}(t)}))= \al_t^{h_i(t)}s\delta_t(\al_s^{h_{ist}(s)})$.

\item[B4.] $Stab_i=\{1,s,tst, stst\}$ (analogue $Stab_i=\{1,t, sts, stst\}$). It holds $i=is, it=its, ist=sts,  itst=itsts$. 
\[
\begin{aligned}
{}[\si_i(s)&\si_{is}(t)\si_{ist}(s) \si_{ists}(t) - \si_{i}(t)\si_{it}(s) \si_{its}(t)\si_{itst}(s) ](f) \\
=&\al_s^{h_i(s)}\delta_s(\al_t^{h_{is}(t)} t(\al_s^{h_{ist}(s)}) t\delta_s(\al_t^{h_{ists}(t)} t(f))) - \al_t^{h_{i}(t)} t(\al_s^{h_{it}(s)}) t\delta_s(\al_t^{h_{its}(t)} t(\al_s^{h_{itst}(s)}) t\delta_s(f)) \\
=& [\al_s^{h_i(s)} \delta_s(\al_t^{h_{is}(t)}t(\al_s^{h_{ist}(s)}) t\delta_s(\al_t^{h_{ists}(t)}))] f \\
+& [t(\al_s^{h_{it}(s)})[s(\al_t^{h_{is}(t)}) st\delta_s(\al_t^{h_{ists}(t)}) - \al_t^{h_{i}(t)}t\delta_s(\al_t^{h_{its}(t)}))] \si_i(s)(f)
\end{aligned}
\]
using $\delta_st\delta_s t = t\delta_s t\delta_s$. 

\end{itemize}
This finishes the investigation of the ten possible cases. We also like to remark that in the example in chapter 5 the case B4 only occurs for 
$Stab_i=\{1,t, sts, stst\}$, i.e. the other stabilizer never occurs. 

\paragraph{(C) $ststst=tststs\colon $ }  $\langle s,t\rangle \cong D_6$, 
\[
\begin{aligned}
t(\al_s)=\al_s+\al_t, &\quad  st(\al_s)= 2\al_s+\al_t, \quad  tst(\al_s)=st(\al_s), \\
s(\al_t)=3\al_s +\al_t, &\quad  ts(\al_t)=3\al_s+2\al_t,\quad  sts(\al_t)= ts(\al_t).
\end{aligned}
\]
It holds 
\[
\begin{aligned}
h_{itstst}(s) =h_i(s) ,& \quad h_{itsts} (t) = h_{is}(t), \quad h_{itst}(s)= h_{ist}(s) \\
h_{its}(t)= h_{ists}(t) , & \quad h_{it}(s) =h_{istst}(s), \quad h_i(t) =h_{iststs}(t).   
\end{aligned}
\]
this implies 
\[
\begin{aligned}
\al_s^{h_{i}(s)}& s(\al_t^{h_{is}(t)}) st (\al_s^{h_{ist}(s)}) sts (\al_t^{h_{ists}(t)}) stst(\al_s^{h_{istst}(s)}) ststs(\al_t^{h_{iststs}(t)}) =\\
\al_t^{h_{i}(t)}& t(\al_s^{h_{it}(s)}) ts(\al_t^{h_{its}(t)}) tst(\al_s^{h_{itst}(s)}) tsts(\al_t^{h_{itsts}(t)}) tstst(\al_s^{h_{itstst}(s)})  
\end{aligned}
\]
Now, $D_6$ has $13$ subgroups. In the following cases the above equality directly implies that $\si_i(ststst)-\si_i(tststs)=0$: 
\[ Stab_i = \{1\},\;\; Stab_i= \{1,tst\},\;\; Stab_i=\{1,sts\},\;\; Stab_i=\{1,ststst\}, \;\; Stab_i= \langle st\rangle =\langle ts\rangle  \]
\begin{itemize}
\item[C1.] $Stab_i=\langle s,t\rangle $. By assumption we have $h_i(s)=0=h_i(t)$ in this case, therefore 
\[
\begin{aligned}
\si_i(s)&\si_i(t)\si_i(s)\si_i(t) \si_i(s) \si_i(t) - \si_i(t)\si_i(s)\si_i(t)\si_i(s)\si_i(t) \si_i(s)= \\
 \delta_s&\delta_t\delta_s\delta_t\delta_s\delta_t - \delta_t\delta_s\delta_t\delta_s\delta_t\delta_s =0
\end{aligned}
\]
because that is known for the divided difference operators, cp \cite{De}. 

\item[C2.] $Stab_i=\{1,s\}$ (analogue $Stab_i=\{1,t\}$). Then, $is=i, itstst=itststs$. 
\[
\begin{aligned}
{}[\si_{i}(s)& \si_{is}(t)\si_{ist}(s)\si_{ists}(t) \si_{istst}(s)\si_{iststs}(t) - \si_{i}(t)\si_{it}(s)\si_{its}(t) \si_{itst}(s)\si_{itsts}(t)\si_{itstst}(s)](f) \\
=&\al_s^{h_{i}(s)}\delta_s(\al_t^{h_{is}(t)} t(\al_s^{h_{ist}(s)}) ts(\al_t^{h_{ists}(t)}) tst (\al_s^{h_{istst}(s)}) tsts(\al_t^{h_{iststs}(t)}) tstst(f)) \\
 &- \al_t^{h_{i}(t)} t(\al_s^{h_{it}(s)}) ts(\al_t^{h_{its}(t)}) tst(\al_s^{h_{itst}(s)}) tsts(\al_t^{h_{itsts}(t)})tstst(\al_s^{h_{itstst}(s)}) tstst\delta_s(f) \\
=&\al_s^{h_{i}(s)}\delta_s(\al_t^{h_{is}(t)}t(\al_s^{h_{ist}(s)}) ts(\al_t^{h_{ists}(t)}) tst(\al_s^{h_{istst}(s)}) tsts(\al_t^{h_{iststs}(t)}) ) tstst(f)\\
=&0
\end{aligned}
\]
using $\delta_s tstst = tstst \delta_s$ and $\delta_s(\al_t^{h_{is}(t)} t(\al_s^{h_{ist}(s)}) ts(\al_t^{h_{ists}(t)}) tst (\al_s^{h_{istst}(s)}) tsts(\al_t^{h_{iststs}(t)}))=0$ because 
\[
\begin{aligned}
s(\al_t^{h_{is}(t)}& t(\al_s^{h_{ist}(s)}) ts(\al_t^{h_{ists}(t)}) tst (\al_s^{h_{istst}(s)}) tsts(\al_t^{h_{iststs}(t)}) )\\
=& \al_t^{h_{is}(t)} t(\al_s^{h_{ist}(s)}) ts(\al_t^{h_{ists}(t)}) tst (\al_s^{h_{istst}(s)}) tsts(\al_t^{h_{iststs}(t)}).
\end{aligned}
\] 

\item[C3.] $Stab_i=\{1,tstst\}$ (analogue $Stab_i=\{1,ststs\}$). Then $its=itst, ists=istst$. 
\[
\begin{aligned}
{}[\si_{i}(s)& \si_{is}(t)\si_{ist}(s)\si_{ists}(t) \si_{istst}(s)\si_{iststs}(t) - \si_{i}(t)\si_{it}(s)\si_{its}(t) \si_{itst}(s)\si_{itsts}(t)\si_{itstst}(s)](f) \\
=&\al_s^{h_{i}(s)}s(\al_t^{h_{is}(t)}) st(\al_s^{h_{ist}(s)}) sts(\al_t^{h_{ists}(t)}) sts \delta_t(\al_s^{h_{istst}(s)} s(\al_t^{h_{iststs}(t)}) st(f)) \\
&-\al_t^{h_{i}(t)} t(\al_s^{h_{it}(s)}) ts(\al_t^{h_{its}(t)}) ts\delta_t(\al_s^{h_{itst}(s)}s(\al_t^{h_{itsts}(t)}) st(\al_s^{h_{itstst}(s)}) sts(f)) \\
=&\al_s^{h_{i}(s)}s(\al_t^{h_{is}(t)}) st(\al_s^{h_{ist}(s)}) sts(\al_t^{h_{ists}(t)}) sts \delta_t(\al_s^{h_{istst}(s)} s(\al_t^{h_{iststs}(t)})) s(f) \\
&+\al_s^{h_{i}(s)}s(\al_t^{h_{is}(t)}) st(\al_s^{h_{ist}(s)}) sts(\al_t^{h_{ists}(t)}) stst(\al_s^{h_{istst}(s)}) ststs(\al_t^{h_{iststs}(t)}) sts\delta_t st(f)\\
&-\al_t^{h_{i}(t)} t(\al_s^{h_{it}(s)}) ts(\al_t^{h_{its}(t)}) ts\delta_t(\al_s^{h_{itst}(s)}s(\al_t^{h_{itsts}(t)}) st(\al_s^{h_{itstst}(s)})) s(f)\\
&-\al_t^{h_{i}(t)} t(\al_s^{h_{it}(s)}) ts(\al_t^{h_{its}(t)}) tst(\al_s^{h_{itst}(s)})tsts(\al_t^{h_{itsts}(t)}) tstst(\al_s^{h_{itstst}(s)}) ts\delta_t sts(f)\\
=& [s(\al_t^{h_{is}(t)}) st(\al_s^{h_{ist}(s)}) sts(\al_t^{h_{ists}(t)}) sts \delta_t(\al_s^{h_{istst}(s)} s(\al_t^{h_{iststs}(t)})) \\
&- \al_t^{h_{i}(t)} t(\al_s^{h_{it}(s)}) ts(\al_t^{h_{its}(t)}) ts\delta_t(\al_s^{h_{itst}(s)}s(\al_t^{h_{itsts}(t)}))] \si_i(s)(f)
\end{aligned}
\]
using $ts\delta_t sts = sts \delta_t st$. 

\item[C4.] $Stab_i=\{1,s, tstst, ststst\}$ ( analogue $Stab_i=\{1,t,ststs, ststst\}$). Then $is=i, itst=its$. Observe, in this case 
\[ h_i(t)=h_{it}(t), \text{ and } h_{it}(s) = h_{its}(s)\]
and it holds 
\[
\begin{aligned}
{}[\si_{i}(s)& \si_{is}(t)\si_{ist}(s)\si_{ists}(t) \si_{istst}(s)\si_{iststs}(t) - \si_{i}(t)\si_{it}(s)\si_{its}(t) \si_{itst}(s)\si_i(tsts)\si_{itstst}(s)](f) \\
=&\al_s^{h_{i}(s)}\delta_s(\al_t^{h_{is}(t)}t(\al_s^{h_{ist}(s)}) ts(\al_t^{h_{ists}(t)}) ts\delta_t(\al_s^{h_{istst}(s)} s(\al_t^{h_{iststs}(t)}) st(f))) \\
&- \al_t^{h_i(t)} t(\al_s^{h_{it}(s)}) ts(\al_t^{h_{its}(t)}) ts\delta_t(\al_s^{h_{itst}(s)}s(\al_t^{h_{itsts}(t)})st(\al_s^{h_{itstst}(s)}) st\delta_s(f)) \\
=& \al_s^{h_{i}(s)}\delta_s(\al_t^{h_{is}(t)}t(\al_s^{h_{ist}(s)}) ts(\al_t^{h_{ists}(t)}) ts\delta_t(\al_s^{h_{istst}(s)} s(\al_t^{h_{iststs}(t)})))\cdot f\\
&+ \al_s^{h_{i}(s)}s(\al_t^{h_{is}(t)})st(\al_s^{h_{ist}(s)}) sts(\al_t^{h_{ists}(t)}) sts\delta_t(\al_s^{h_{istst}(s)} s(\al_t^{h_{iststs}(t)})) \delta_s(f)\\
&+\al_s^{h_{i}(s)}\delta_s(\al_t^{h_{is}(t)}t(\al_s^{h_{ist}(s)}) ts(\al_t^{h_{ists}(t)}) tst(\al_s^{h_{istst}(s)}) tsts(\al_t^{h_{iststs}(t)})) ts\delta_t st(f)))\\
&+\al_s^{h_{i}(s)} s(\al_t^{h_{is}(t)})st(\al_s^{h_{ist}(s)}) sts(\al_t^{h_{ists}(t)}) stst(\al_s^{h_{istst}(s)}) ststs(\al_t^{h_{iststs}(t)}) \delta_s ts \delta_t st(f)\\
&- \al_t^{h_i(t)} t(\al_s^{h_{it}(s)}) ts(\al_t^{h_{its}(t)}) ts\delta_t(\al_s^{h_{itst}(s)}s(\al_t^{h_{itsts}(t)})st(\al_s^{h_{itstst}(s)})) \delta_s(f)) \\
&-\al_t^{h_i(t)} t(\al_s^{h_{it}(s)}) ts(\al_t^{h_{its}(t)}) tst(\al_s^{h_{itst}(s)})tsts(\al_t^{h_{itsts}(t)})tstst(\al_s^{h_{itstst}(s)}) ts\delta_t st\delta_s(f)) \\
&= [\al_s^{h_{i}(s)}\delta_s(\al_t^{h_{is}(t)}t(\al_s^{h_{ist}(s)}) ts(\al_t^{h_{ists}(t)}) ts\delta_t(\al_s^{h_{istst}(s)} s(\al_t^{h_{iststs}(t)})))] f\\
&+[ s(\al_t^{h_{is}(t)})st(\al_s^{h_{ist}(s)}) sts(\al_t^{h_{ists}(t)}) sts\delta_t(\al_s^{h_{istst}(s)} s(\al_t^{h_{iststs}(t)})) \\
&\; - \al_t^{h_i(t)} t(\al_s^{h_{it}(s)}) ts(\al_t^{h_{its}(t)}) ts\delta_t(\al_s^{h_{itst}(s)}s(\al_t^{h_{itsts}(t)}))] \si_i(s) (f) 
\end{aligned}
\]
using $\delta_s ts \delta_t st = ts\delta_t st \delta_s$. 

\item[C5.] $Stab_i=\{ 1,sts,tst, ststst\}$. Then $is = ist, it=its$. Observe, in this case 
\[ h_i(s)= h_{is}(s), \text{ and } h_i(t)= h_{it}(t) \]
and it holds 
\[
\begin{aligned}
{}[\si_{i}(s)& \si_{is}(t)\si_{ist}(s)\si_{ists}(t) \si_{istst}(s)\si_{iststs}(t) - \si_{i}(t)\si_{it}(s)\si_{its}(t) \si_{itst}(s)\si_{itsts}(t)\si_{itstst}(s)](f) \\
=&\al_s^{h_{i}(s)} s(\al_t^{h_{is}(t)}) s\delta_t(\al_s^{h_{ist}(s)} s(\al_t^{h_{ists}(t)}) st(\al_s^{h_{istst}(s)})st\delta_s(\al_t^{h_{iststs}(t)} t(f))) \\
&- \al_t^{h_i(t)} t(\al_s^{h_{it}(s)}) t\delta_s(\al_t^{h_{its}(t)}t(\al_s^{h_{itst}(s)}) ts(\al_t^{h_{itsts}(t)}) ts\delta_t(\al_s^{h_{itstst}(s)} s(f))) \\
=&\al_s^{h_{i}(s)} s(\al_t^{h_{is}(t)}) s\delta_t(\al_s^{h_{ist}(s)} s(\al_t^{h_{ists}(t)}) st(\al_s^{h_{istst}(s)})st\delta_s(\al_t^{h_{iststs}(t)}))\cdot  f \\
&+\al_s^{h_{i}(s)} s(\al_t^{h_{is}(t)}) st(\al_s^{h_{ist}(s)}) sts(\al_t^{h_{ists}(t)}) stst(\al_s^{h_{istst}(s)})stst\delta_s(\al_t^{h_{iststs}(t)}) s\delta_t s(f) \\
&+\al_s^{h_{i}(s)} s(\al_t^{h_{is}(t)}) s\delta_t(\al_s^{h_{ist}(s)} s(\al_t^{h_{ists}(t)}) st(\al_s^{h_{istst}(s)})sts(\al_t^{h_{iststs}(t)})) t\delta_s t(f) \\
&+\al_s^{h_{i}(s)} s(\al_t^{h_{is}(t)}) st(\al_s^{h_{ist}(s)}) sts(\al_t^{h_{ists}(t)}) stst(\al_s^{h_{istst}(s)})ststs(\al_t^{h_{iststs}(t)}) s\delta_t st \delta_s t(f) \\
&- \al_t^{h_i(t)} t(\al_s^{h_{it}(s)}) t\delta_s(\al_t^{h_{its}(t)}t(\al_s^{h_{itst}(s)}) ts(\al_t^{h_{itsts}(t)}) ts\delta_t(\al_s^{h_{itstst}(s)})) \cdot f \\
&-\al_t^{h_i(t)} t(\al_s^{h_{it}(s)}) ts(\al_t^{h_{its}(t)}) tst(\al_s^{h_{itst}(s)}) tsts(\al_t^{h_{itsts}(t)}) tsts\delta_t(\al_s^{h_{itstst}(s)}) t\delta_st(f) \\
&-\al_t^{h_i(t)} t(\al_s^{h_{it}(s)}) t\delta_s(\al_t^{h_{its}(t)}t(\al_s^{h_{itst}(s)}) ts(\al_t^{h_{itsts}(t)}) tst(\al_s^{h_{itstst}(s)})) s\delta_t s(f) \\
&-\al_t^{h_i(t)} t(\al_s^{h_{it}(s)}) ts(\al_t^{h_{its}(t)}tst(\al_s^{h_{itst}(s)}) tsts(\al_t^{h_{itsts}(t)}) tstst(\al_s^{h_{itstst}(s)}) t\delta_s ts \delta_t s(f) \\
=& P_e f + [ \al_s^{h_{i}(s)} s(\al_t^{h_{is}(t)})st(\al_s^{h_{ist}(s)}) sts(\al_t^{h_{ists}(t)}) stst(\al_s^{h_{istst}(s)})stst\delta_s(\al_t^{h_{iststs}(t)})  \\
&- \al_t^{h_i(t)} t(\al_s^{h_{it}(s)}) t\delta_s(\al_t^{h_{its}(t)}t(\al_s^{h_{itst}(s)}) ts(\al_t^{h_{itsts}(t)}) tst(\al_s^{h_{itstst}(s)}))] s\delta_t s (f) \\
&+[\al_s^{h_{i}(s)} s(\al_t^{h_{is}(t)}) s\delta_t(\al_s^{h_{ist}(s)} s(\al_t^{h_{ists}(t)}) st(\al_s^{h_{istst}(s)})sts(\al_t^{h_{iststs}(t)}))\\
&- \al_t^{h_i(t)} t(\al_s^{h_{it}(s)}) ts(\al_t^{h_{its}(t)}) tst(\al_s^{h_{itst}(s)}) tsts(\al_t^{h_{itsts}(t)}) tsts\delta_t(\al_s^{h_{itstst}(s)})] t\delta_s t (f)
\end{aligned}
\]
using $s\delta_t st \delta_s t = t\delta_s ts \delta_t s $ where $P_e$ is a polynomial in $\al_t, \al_s$. Then we look at 
\[ 
\begin{aligned}
\si_{i}(s)\si_{is}(t) \si_{ist} (s) (f) &= \al_s^{h_{i}(s)} s(\al_t^{h_{is}(t)}) s\delta_t(\al_s^{h_{ist}(s)}) \cdot f + \al_s^{h_i(s)} s(\al_t^{h_{is}(t)}) st(\al_s^{h_{ist}(s)}) s\delta_t s(f) \\
\si_{i}(t)\si_{it}(s) \si_{its} (t) (f) &= \al_t^{h_{i}(t)} t(\al_s^{h_{it}(s)}) t\delta_s(\al_t^{h_{its}(t)}) \cdot f + \al_t^{h_i(t)} t(\al_s^{h_{it}(s)}) ts(\al_t^{h_{its}(t)}) t\delta_s t(f)
\end{aligned}
\]  
and we observe for the coefficient in front of $s\delta_t s(f)$ that it is divisible by \\
$\al_s^{h_i(s)} s(\al_t^{h_{is}(t)}) st(\al_s^{h_{ist}(s)})$ and the one in front of $t\delta_s t$ is divisible by $\al_t^{h_i(t)} t(\al_s^{h_{it}(s)}) ts(\al_t^{h_{its}(t)})$. Observe $h_{ist}(s)= h_i(s), h_{its}(t)=h_i(t)$. Use the following simplifictaions 
\[
\begin{aligned}
t\delta_s(\al_t^{h_{its}(t)}&t(\al_s^{h_{itst}(s)}) ts(\al_t^{h_{itsts}(t)}) tst(\al_s^{h_{itstst}(s)})) \\
=& t\delta_s(\al_t^{h_i(t)}t(\al_s^{h_i(s)}) ts(\al_t^{h_{is}(t)}) tst(\al_s^{h_i(s)})) \\
=&s(\al_t^{h_{is}(t)}) t\delta_s(\al_t^{h_i(t)} t(\al_s^{h_i(s)}) st(\al_s^{h_i(s)})) \\
=&s(\al_t^{h_{is}(t)})[t\delta_s(\al_t^{h_i(t)} t(\al_s^{h_i(s)})) tst(\al_s^{h_i(s)})+ ts(\al_t^{h_i(t)}) tst(\al_s^{h_i(s)}) t\delta_s st(\al_s^{h_i(s)})]\\
=&s(\al_t^{h_{is}(t)})st(\al_s^{h_i(s)})[\al_s^{h_i(s)}t\delta_s(\al_t^{h_i(t)})+ ts(\al_t^{h_i(t)}) t\delta_s t(\al_s^{h_i(s)}) \\
&\quad - ts(\al_t^{h_i(t)}) t\delta_s t(\al_s^{h_i(s)})]\\
=& \al_s^{h_i(s)}s(\al_t^{h_{is}(t)})st(\al_s^{h_i(s)})t\delta_s(\al_t^{h_i(t)})
\end{aligned}
\]
and analogously 
\[
\begin{aligned}
s\delta_t(\al_s^{h_{ist}(s)}&s(\al_t^{h_{ists}(t)}) st(\al_s^{h_{istst}(s)}) sts(\al_t^{h_{iststs}(t)}))\\
 =&s\delta_t(\al_s^{h_i(s)}s(\al_t^{h_i(t)}) st(\al_s^{h_{it}(s)}) sts(\al_t^{h_i(t)}))\\ =&\al_t^{h_i(t)}t(\al_s^{h_{it}(s)})ts(\al_t^{h_i(t)})s\delta_t(\al_s^{h_i(s)})
\end{aligned}
\]
Then a simple substitution gives that the difference above is of the form 
\[Q_e f + Q_{sts} \si_{i}(s)\si_{is}(t) \si_{ist} (s)(f) + Q_{tst}\si_{i}(t)\si_{it}(s) \si_{its} (t)\]
for some polynomials $Q_e, Q_{sts}, Q_{tst}$ in $\al_s, \al_t$.
\end{itemize}

\end{itemize}
Now, let $A$ be the algebra given by generator $\widetilde{1_i}, \widetilde{z_i(t)}, \widetilde{\si_i(s)}$ subject to relations (1)-(5). 
Then, by the straightening rule and the braid relation it holds that if $w=s_1\cdots s_k =t_1\cdots t_k$ are two reduced expressions then 
\[ \widetilde{\si (t_1\cdots t_k)}\in \sum_{v\leq s_1\cdots s_k \text{ reduced subword }} \mcE * \widetilde{\si (v)}. \] 
Therefore, once we have fixed one (any) reduced expression for each for $w \in \W$, it holds 
\[ A=\sum_{w\in \W} \mcE * \widetilde{\si(w)}.\]
Since the generators of $\mcZ_G$ fulfill the relations (1)-(5), we have a surjective algebra homomorphism 
\[ A\to \mcZ_G \]
mapping $\widetilde{1_i}\mapsto 1_i, \widetilde{z_i(t)}\mapsto z_i(t), \widetilde{\si_i(s)}\mapsto \si_i(s)$. 
Since $\mcZ_G= \bigoplus_{w\in \W} \mcE* \si (w)$ and the map is by definition $\mcE$-linear 
it follows that $A= \bigoplus_{w\in \W} \mcE * \widetilde{\si(w)}$ and the map is an isomorphism. 
\hfill $\Box$
\addcontentsline{toc}{section}{\textbf{References} \hfill}

\bibliographystyle{alphadin}
\bibliography{GenQuiverHecke}

 \end{document}